\documentclass{article}

\usepackage{amssymb}
\usepackage{amsmath}
\usepackage{amsthm}
\usepackage{mathtools}
\usepackage{graphicx}
\usepackage{epstopdf}
\usepackage{caption}
\usepackage{multirow}
\usepackage{booktabs}
\usepackage[flushleft]{threeparttable}
\usepackage{adjustbox}
\usepackage{array}
\usepackage{bm}
\usepackage[svgnames]{xcolor}
\usepackage{url}
\usepackage[colorlinks=true, linkcolor=Brown, citecolor=DarkCyan, urlcolor=Olive]{hyperref}
\usepackage{array}
\newcolumntype{C}{>{$}c<{$}}
\usepackage{appendix}
\usepackage[square, numbers]{natbib}

\usepackage{xspace}
\usepackage{siunitx}
\usepackage{newfloat}
\usepackage{chngcntr}
\usepackage{mathtools}
\usepackage{breqn}

\DeclareFloatingEnvironment{snippet}

\usepackage{xcolor,colortbl}
\definecolor{Gray}{gray}{0.85}
\definecolor{LightGray}{gray}{0.95}
\newcolumntype{a}{>{\columncolor{Gray}}c}
\newcolumntype{b}{>{\columncolor{LightGray}}c}

\RequirePackage{bm}
\RequirePackage{url}
\newcommand{\vectorFont}[1]{\boldsymbol{#1}}

\newcommand{\matrixFont}[1]{\mathbf{#1}}

\newcommand{\dd}[2]{\ensuremath{\frac{\partial #1}{\partial #2}}\xspace}

\newcommand{\scalar}[2]{\ensuremath{\zeta^{#1}_{#2}}\xspace}
\newcommand{\vect}[1][]{\ensuremath{\vectorFont{\iota}_{#1}}\xspace}

\newcommand{\temperature}[2]{\ensuremath{T^{#1}_{#2}}\xspace}

\newcommand{\internalEnergy}[2]{\ensuremath{u^{#1}_{#2}}\xspace}
\newcommand{\enthalpy}[2]{\ensuremath{h^{#1}_{#2}}\xspace}

\newcommand{\displacement}[2]{\ensuremath{\vectorFont{u}^{#1}_{#2}}\xspace}
\newcommand{\pressure}[2]{\ensuremath{p^{#1}_{#2}}\xspace}
\newcommand{\tangentialComponent}{\ensuremath{\parallel}\xspace}
\newcommand{\normalComponent}{\ensuremath{\perp}\xspace}
\newcommand{\traction}[2]{\ensuremath{\vectorFont{\lambda}^{#1}_{#2}}\xspace}
\newcommand{\normalTraction}{\ensuremath{\lambda_{\normalComponent}}\xspace}
\newcommand{\fluidVelocity}[2]{\ensuremath{\vectorFont{v}^{#1}_{#2}}\xspace}
\newcommand{\fluidFlux}[2]{\ensuremath{\vectorFont{v}^{#1}_{#2}}\xspace}

\newcommand{\density}[2]{\ensuremath{\rho^{#1}_{#2}}\xspace}

\newcommand{\stress}[2]{\ensuremath{\sigma^{#1}_{#2}}\xspace}

\newcommand{\enthalpyFlux}[2]{\ensuremath{\vectorFont{w}^{#1}_{#2}}\xspace}
\newcommand{\heatFlux}[2]{\ensuremath{\vectorFont{q}^{#1}_{#2}}\xspace}

\newcommand{\interfaceEnthalpyFlux}[2]{\ensuremath{w^{#1}_{#2}}\xspace}
\newcommand{\interfaceHeatFlux}[2]{\ensuremath{q^{#1}_{#2}}\xspace}
\newcommand{\interfaceFluidFlux}[2]{\ensuremath{v^{#1}_{#2}}\xspace}

\newcommand{\subdomain}[2]{\ensuremath{\Omega^{#1}_{#2}}\xspace}
\newcommand{\boundary}[2]{\ensuremath{\partial_{#2}\Omega_{#1}}\xspace}
\newcommand{\interface}[2]{\ensuremath{\Gamma^{#1}_{#2}}\xspace}

\newcommand{\higherSet}[1][]{\ensuremath{\hat{S}_{#1}}\xspace}

\newcommand{\aperture}[2]{\ensuremath{a^{#1}_{#2}}\xspace}
\newcommand{\specificVolume}[2]{\ensuremath{\mathcal{V}^{#1}_{#2}}\xspace}
\newcommand{\normal}[2]{\ensuremath{\vectorFont{n}^{#1}_{#2}}\xspace}

\newcommand{\gap}[2]{\ensuremath{g^{#1}_{#2}}\xspace}

\newcommand{\massSource}[2]{\ensuremath{\psi^{#1}_{#2}}\xspace}
\newcommand{\energySource}[2]{\ensuremath{\chi^{#1}_{#2}}\xspace}

\newcommand{\bodyForce}[2]{\ensuremath{\vectorFont{F}^{#1}_{#2}}\xspace}

\newcommand{\bulkModulus}[1]{\ensuremath{K_{#1}}\xspace}
\newcommand{\shearModulus}{\ensuremath{G}\xspace}

\newcommand{\gravityVector}[1][]{\ensuremath{\vectorFont{g}}\xspace}%\vectorFont{g}_{#1}}\xspace}

\newcommand{\thermalExpansion}[2]{\ensuremath{\beta^{#1}_{#2}}\xspace}
\newcommand{\thermalConductivity}[2]{\ensuremath{\kappa^{#1}_{#2}}\xspace}
\newcommand{\heatCapacity}[2]{\ensuremath{c^{#1}_{#2}}\xspace}

\newcommand{\biotAlpha}{\ensuremath{\alpha}\xspace}
\newcommand{\porosity}[2]{\ensuremath{\phi^{#1}_{#2}}\xspace}
\newcommand{\compressibility}[2]{\ensuremath{\gamma^{#1}_{#2}}\xspace}
\newcommand{\permeability}[2]{\ensuremath{\mathcal{K}^{#1}_{#2}}\xspace}

\newcommand{\viscosity}[2]{\ensuremath{\eta^{#1}_{#2}}\xspace}
\newcommand{\frictionCoefficient}{\ensuremath{F}\xspace}
\newcommand{\frictionBound}{\ensuremath{b}\xspace}

\newcommand{\dilationAngle}{\ensuremath{\psi}\xspace}
\newcommand{\maximumFractureClosure}{\ensuremath{\Delta u_{max}}\xspace}
\newcommand{\normalStiffness}{\ensuremath{K_n}\xspace}
\newcommand{\incrementTangentialDisplacement}[1][]{\ensuremath{\jump{\dot{\displacement{}{}}^{#1}}_{\tangentialComponent}}\xspace}
\newcommand{\identity}[1][]{\ensuremath{\matrixFont{I}_{#1}}\xspace}
\newcommand{\norm}[1]{\ensuremath{||#1||}\xspace}
\newcommand{\jump}[1]{\ensuremath{[\![#1]\!]}\xspace}
\newcommand{\effective}[1]{\ensuremath{\{#1}\}\xspace}

\newcommand{\projectFromInterface}[2]{\ensuremath{\Xi ^{#1}_{#2}}\xspace}
\newcommand{\projectToInterface}[2]{\ensuremath{\Pi^{#1}_{#2}}\xspace}
\newcommand{\trace}[1]{\ensuremath{\text{tr}({#1}) }\xspace}
\newcommand{\plusSide}{\ensuremath{j}\xspace}
\newcommand{\minusSide}{\ensuremath{k}\xspace}

\DeclareSIUnit\year{yr}

\begin{document}
% Title
\title{\vspace{-1.2cm} Flexible and rigorous numerical modelling of multiphysics processes in fractured porous media using PorePy}

% Author's information
\author{Ivar Stefansson\footnote{Corresponding author e-mail: \href{mailto:ivar.stefansson@uib.no}{ivar.stefansson@uib.no.}} $^{,}$\footnote{Center for Modeling of Coupled Subsurface Dynamics, University of Bergen, P.O. Box 7800, N-5020 Bergen, Norway.} \and Jhabriel Varela$^{\dagger}$ \and Eirik Keilegavlen$^{\dagger}$ \and Inga Berre$^{\dagger}$}

\date{\today}

% Insert title
\maketitle

% Abstract
\begin{abstract}
\thispagestyle{empty}
    Multiphysics processes in fractured porous media is a research field of importance for several subsurface applications and has received considerable attention over the last decade.
    The dynamics are characterised by strong couplings between processes as well as interaction between the processes and the structure of the fractured medium itself.
    The rich range of behavior calls for explorative mathematical modelling, such as experimentation with constitutive laws and novel coupling concepts between physical processes.
    Moreover, efficient simulations of the strong couplings between multiphysics processes and geological structures require the development of tailored numerical methods.

    We present a modelling framework and its implementation in the open-source simulation toolbox \texttt{PorePy}, which is designed for rapid prototyping of multiphysics processes in fractured porous media.
    \texttt{PorePy} uses a mixed-dimensional representation of the fracture geometry and generally applies fully implicit couplings between processes.
    The code design follows the paradigms of modularity and differentiable programming, which together allow for extreme flexibility in experimentation with governing equations with minimal changes to the code base.
    The code integrity is supported by a multilevel testing framework ensuring the reliability of the code.

    We present our modelling framework within a context of thermo-poroelasticity in deformable fractured porous media, illustrating the close relation between the governing equations and the source code.
    We furthermore discuss the design of the testing framework and present simulations showcasing the extendibility of \texttt{PorePy}, as well as the type of results that can be produced by mixed-dimensional simulation tools.\vspace{3mm}  \vspace{3mm}\\
\end{abstract}

\section*{Highlights}

\begin{itemize}
    \item PorePy simulation toolbox for rapid prototyping of multiphysics processes in fractured porous media
    \item Discretized model adheres closely to governing equations
    \item Flexibility ensured by modular implementation and differentiable programming
    \item Integrated testing framework
\end{itemize}
\thispagestyle{empty}
\tableofcontents

\newpage

\section{Introduction}
The context for this work is the modelling and simulation of multiphysics processes in fractured porous media,
with applications including extraction of geothermal energy and hydrocarbons, wastewater disposal, $\mathrm{CO_2}$ storage and storage of chemical and thermal energy.
Noting that several aspects of the discussion apply to other settings as well, the primary focus in this paper is on thermo-hydro-mechanical processes.
These processes not only involve complex and coupled dynamics.
Observing and monitoring them is inherently difficult due to their subsurface situation in heterogeneous and fractured rock formations.
This motivates the use of mathematical and numerical modelling, thus creating demand for simulation toolboxes for the above-mentioned class of problems.
Due to tight coupling between different processes and severe structural impact of fractures, standard simulators have limitations in their applicability, leading to active development of research codes.

A research simulation tool for this context may serve at least two purposes.
First, it can facilitate simulation technology research into tailored numerical methods.
Second, it can accommodate modelling studies enhancing our understanding of the processes at play.
Furthermore, we single out two important properties which characterise a toolbox which is fit for those two purposes.
First, rigorous adherence to the mathematical model is required to ensure that the governing physical and constitutive laws are fulfilled up to discretization error.
This should be understood broadly and includes solving the full system of governing equations without decoupling simplifications, employing sound discretization schemes and thorough and structured code testing.
Second, a workable tool must be flexible to accommodate relevant adaptation and extension.
This naturally includes the reasonable requirement for any simulation tool of convenient problem specification through geometry, material parameters, boundary conditions and the like.
However, for modelling research it extends to modifications of governing equations.
Similarly, research into the simulation technology itself requires choosing and experimenting with meshing, discretization schemes, solvers for non-linear and linear equation systems, etc.
Ideally, all this flexibility should be accommodated while minimising code complexity and maintaining user-friendliness.

A number of simulation toolboxes exist for thermo-poromechanical porous media processes, some of which include fracture representation.
A non-exhaustive list includes \texttt{CSMP}~\cite{matthai2001csmp}, \texttt{MRST}~\cite{lie2019mrst_book}, \texttt{TOUGH3}~\cite{jung2017tough3}, \texttt{GEOSX}~\cite{gross2021geosx}, \texttt{FALCON}~\cite{podgorney2021falcon}, \texttt{GOLEM}~\cite{cacace2017flexible}, \texttt{DuMux}~\cite{koch2021dumux}, \texttt{OpenGeoSys}~\cite{ogs:6.4.3}, \texttt{DARTS}~\cite{voskov2017darts} and \texttt{Flow123d}~\cite{bvrezina2016flow123d}.
Benchmark studies by White et al.~\cite{white2018benchmark_suite} and Mindel et al.~\cite{mindel2021benchmark} provide demonstration of
capabilities and relative strengths and weaknesses of several of these toolboxes.

In this paper, we present features of the \texttt{PorePy} toolbox \cite{PorePy}, an open-source research code, written in Python and tailored to the specific needs of simulating multiphysics in fractured porous media.
We focus on its rigorous adherence to the mathematical model and how it accommodates flexibility in model setups and extensions.
\texttt{PorePy}'s main design principles are presented in Section \ref{sec:design:principles},\ while the mathematical model and its implementation is presented in Section \ref{sec:mathematical:models}.\
Numerical solution strategies are presented in Section \ref{sec:numerical}.\
We present a selected suite of tests verifying the code's adherence to a mixed-dimensional mathematical model for compressible flow in Section \ref{sec:testing},\ before demonstrating code versatility in application to coupled processes  in Section \ref{sec:application:examples}. Finally, in Section \ref{sec:conclusion}, we present our concluding remarks.

\section{Design principles\label{sec:design:principles}}
This section discusses the high-level code design in \texttt{PorePy} describing how it promotes the requirements of rigour and flexibility identified in the introduction.
We commence with an overview of the code's core structure and important components,
some of which are expounded on with examples and illustrated using snippets in subsequent sections.

Since fractures have very high aspect ratios, it is natural to represent them as objects one dimension lower than the host domain. Correspondingly, the simulation domain is decomposed into subdomains of dimension successively reduced by one for rock matrix, individual fractures, and fracture intersections. Additionally, each pair of geometrically neighbouring subdomains one dimension apart is connected through an interface. The grids for individual subdomains and interfaces are collected in a graph representing the mixed-dimensional grid.

The implementation relates closely to previously presented models for fractured porous media \cite{Martin2005, boon2018robust, nordbotten2019unified} and mixed-dimensional theory based on exterior calculus \cite{boon2021mixed_dimensional, boon2022mixed}.
This constitutes a rigorous foundation for the formulation of discrete mixed-dimensional models with a close relationship to theoretical results. As such, the implementation is well suited for testing and verification, as demonstrated by Varela et al.\ \cite{varela2022aposteriori}, who tested different local mass-conservative schemes using mixed-dimensional \textit{a posteriori} error estimates.

We define variables, discretizations and arrays representing material parameters and the like on sets of subdomain or interface grids.
All of these elements are compatible with \texttt{PorePy}'s automatic differentiation (AD) framework,
which allows the construction of higher-order elements such as equations by arithmetic operations following a differentiable programming paradigm. On the highest level, the equations are derived from fundamental conservation principles and kinematic constraints and may be combined to compose a multiphysics system.
These include conservation equations for mass, momentum and energy and kinematic constraints for fracture contact.
The fundamental principles are supplemented by constitutive equations prescribing relationships between primary and secondary variables.

Modelling flexibility is achieved through a modular design, allowing all terms of the governing equations to be defined and replaced independently.
Thus, modelling choices can be made by selecting from the options in \texttt{PorePy}'s library of constitutive laws.
The decomposition into subdomains also facilitates using different laws for different subdomain sets, even within the same dimension.
Analogously, individual steps of the solution strategy can be tailored for experimentation with simulation technology.
Extension of the code base is similarly achieved by defining the new relationship and combining it with existing code.
\texttt{Python} being a high-level language, such extensions are also straightforward and do not require expert knowledge of \texttt{PorePy}.

Setting up a complete mixed-dimensional multiphysics model is a nontrivial task.
\texttt{PorePy} therefore provides a suite of ready-to-run setups which we refer to as \texttt{Models}, see Table \ref{tab:models} for a list of the currently covered cases.
As will be discussed in the sequel, the \texttt{Models} provide a base for setting up simulations that are convenient in facilitating flexible code reuse and extension and also are thoroughly tested.

\begin{table}
    \centering
    \caption{\texttt{PorePy} contains \texttt{Models} for the following combinations of conservation laws, with the specific physical model arising from the default choices of constitutive laws specified to the right.\label{tab:models}}
    \begin{tabular}{|l|p{4cm}|p{5cm}|}
    \hline
    & Conservation equation(s) & Default physical model \\
    \hline
        \multirow{2}{1.7cm}{single-physics} & Fluid mass & Compressible single-phase flow \\
        \cline{2-3}
        &Momentum & Elasticity with contact mechanics \\
        \hline
        \multirow{3}{1.7cm}{multi-physics} & Fluid mass and energy  & Compressible single-phase flow \\
        \cline{2-3}
        & Fluid mass and momentum & Poroelasticity with contact mechanics \\
        \cline{2-3}
        &Fluid mass, energy and momentum & Thermo-poroelasticity with contact mechanics \\
        \hline
    \end{tabular}
\end{table}

\section{Mathematical model}
\label{sec:mathematical:models}
In specifying the mathematical model, we present the representation of the mixed-dimensional geometry followed by the equations representing conservation laws and kinematic constraints.
Finally, we close the model by specifying constitutive laws together with initial and boundary conditions.
We illustrate concepts discussed in the previous section using snippets for selected parts of the model.
\subsection{Mixed-dimensional geometry}
\label{sec:mixed-dimensional:geometry}
Reflecting the mixed-dimensional theory for fractured porous media \cite{Martin2005,boon2021mixed_dimensional}, we describe the fractured porous medium as a collection of subdomains \subdomain{}{i} of different dimension $d_i$, with $d_i \in \{0, \ldots, N\}$ and $N \in \{2, 3\}$.
Focusing on the case $N=3$, we represent the porous medium by a 3D subdomain, fractures by 2D subdomains, fracture intersections by 1D subdomains, and intersections of fracture intersections by 0D subdomains.
The width of a dimensionally reduced fracture is characterised by its aperture \aperture{}{i}.
To account for the reduced dimensions of the various subdomains in the full mixed-dimensional setting, we define a
 specific volume \specificVolume{}{i} having dimensions $\si{\meter}^{N-d_i}$ \cite{boon2018robust}.
The relation between \aperture{}{i} and \specificVolume{}{i} will be detailed below.

\begin{figure}
    \centering
    \includegraphics[width=0.38\textwidth]{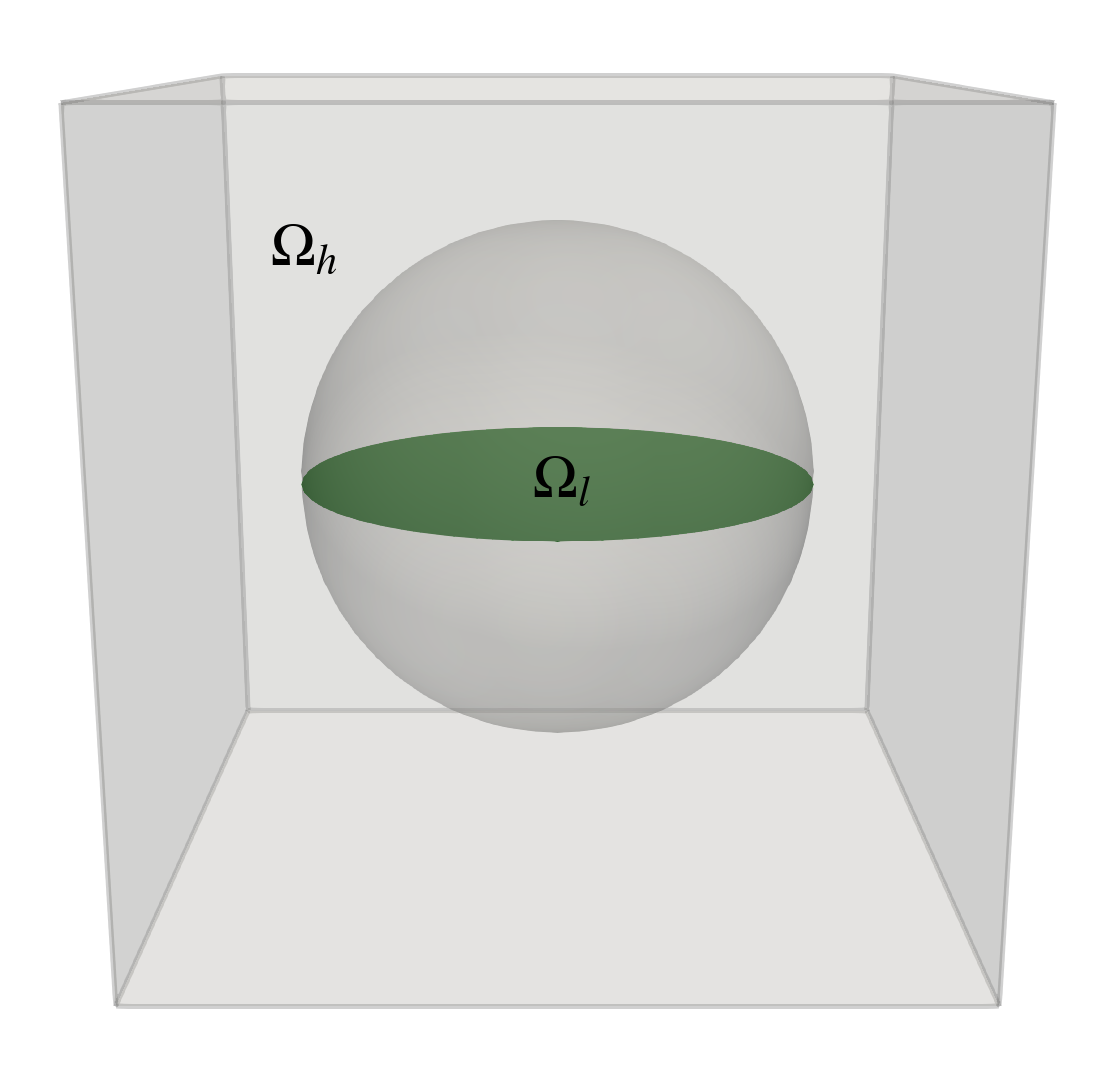}
    \hspace{0.03\textwidth}
    \includegraphics[width=0.28\textwidth]{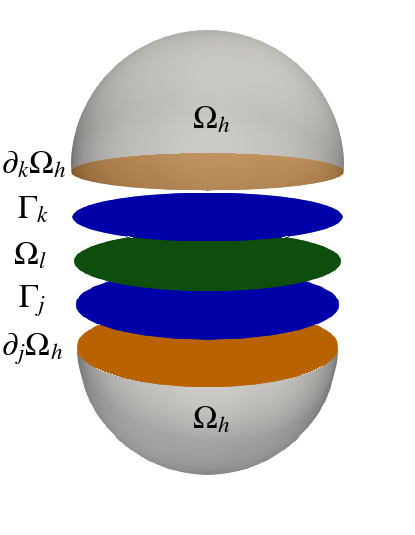}
    \includegraphics[width=0.28\textwidth]{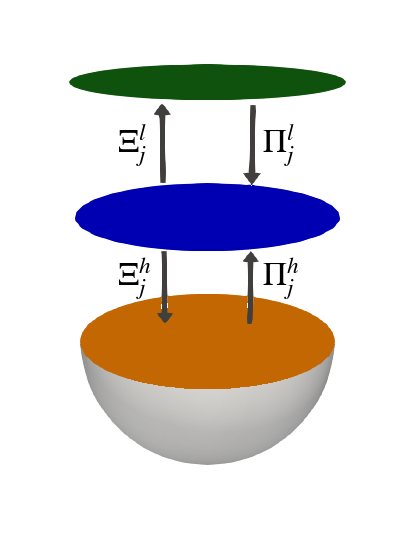}
    \caption{Illustration of a mixed-dimensional geometry.
    To the left, we show the full geometry with a matrix subdomain \subdomain{}{h} and a single circular fracture \subdomain{}{l}.\
    The transparent sphere indicates the area shown in the close-up in the central figure.
    The close-up illustrates the fracture (green), interfaces on either side (blue) and boundaries (orange), all separated for visualisation purposes.
    To the right, we show the projection operators corresponding to the bottom-half coupling between \subdomain{}{h} and \subdomain{}{l} through \interface{}{j}.
    \label{fig:geometry}}
\end{figure}

An interface \interface{}{j} facilitates coupling between each pair of subdomains one dimension apart, see Figure \ref{fig:geometry} for illustration.
We write \boundary{i}{} for the boundary of \subdomain{}{i} and the internal part geometrically coinciding with the interface \interface{}{j} is \boundary{i}{j} $\subseteq \boundary{i}{}$.
We denote the projection of relevant quantities  from subdomain \subdomain{}{i} to interface \interface{}{j} by \projectToInterface{i}{j} and the reverse operation by \projectFromInterface{i}{j}, see the right panel of Figure \ref{fig:geometry}.

We use subscripts to identify quantities associated to subdomains and interfaces and superscripts $f$ and $s$ to denote respectively fluid and solid quantities. However, when context allows, we will suppress subscripts and superscripts in the interest of readability.

\subsection{Conservation laws}
\label{sec:mathematical:models:conservation_laws}
In the following, we present conservation laws for mass, energy and momentum for the relevant subdomains.

The fluid mass conservation equation for a subdomain \subdomain{}{i} of dimension $d_i \in \{0, \ldots, N\}$ reads
\begin{equation}
    \label{eq:mass_conservation}
    \dd{\left(\specificVolume{}{i} \density{f}{i}\porosity{}{i}\right)}{t}
    + \nabla \cdot \left( \specificVolume{}{i}\density{f}{i}\fluidVelocity{}{i} \right)
    - \sum_{j\in\higherSet[i]} \projectFromInterface{i}{j}  \left(\specificVolume{}{j}\density{f}{j}\interfaceFluidFlux{}{j}\right)
    = \massSource{}{i},
\end{equation}
where \density{f}{} is the fluid density, \porosity{}{i} is the porosity, \fluidVelocity{}{i} and \interfaceFluidFlux{}{j} are subdomain and interface volumetric fluid fluxes and \massSource{}{i} is a source or sink of fluid mass.
$\specificVolume{}{j} := \projectToInterface{i}{j} \specificVolume{}{i}$ is the interface specific volume, while
the set \higherSet[i] contains all interfaces to higher-dimensional neighbours of \subdomain{}{i}.\
The second term is void for $d_i = 0$ (since there are no mass fluxes associated with intersection points), whereas the third term is void for $d_i = N$ (as $N$-dimensional subdomains do not have higher-dimensional neighbours).
The \texttt{PorePy} mass balance equation method is shown in Snippet~\ref{code:mass_balance_equation}.
Notice how each term is specified in a separate method, thus facilitating tailoring of any one term with minimal changes to the code.

\begin{snippet}
\caption{The mass balance equation  \eqref{eq:mass_conservation} is assembled by passing all terms to a generic balance equation method. Here, \texttt{pp} refers to the \texttt{PorePy} library and the meaning of the term \texttt{pp.ad.Operator} will be explained in Section \ref{sec:numerical}.}
\fbox{\includegraphics[width=0.975\textwidth]{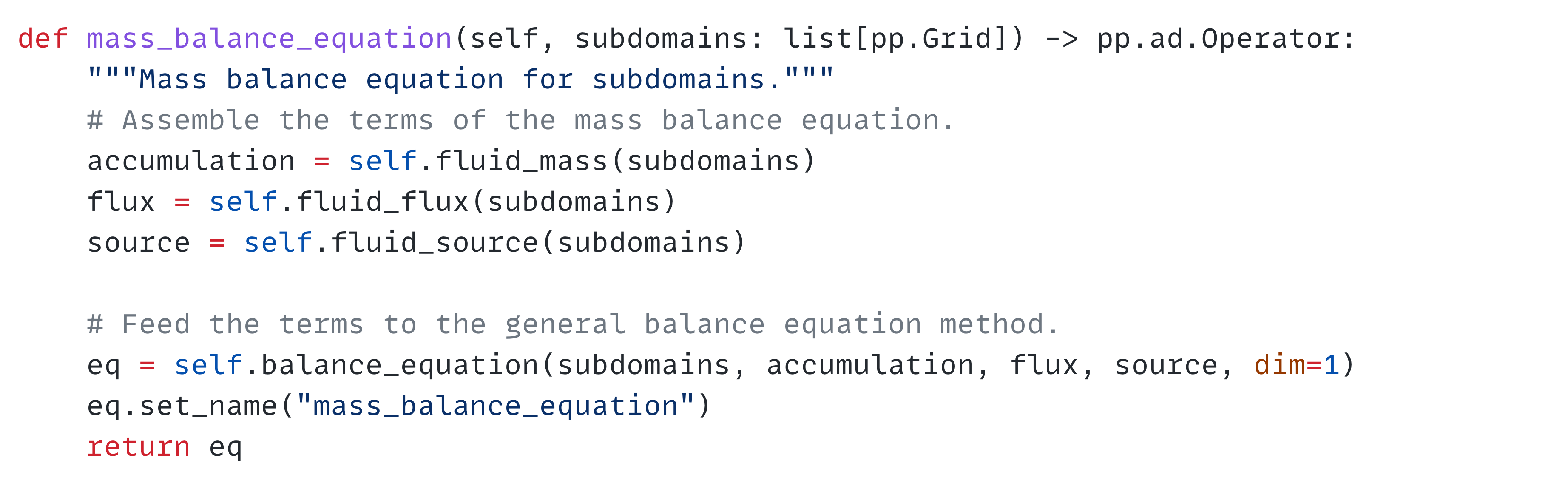}} \label{code:mass_balance_equation}
\end{snippet}

Still considering $d_i \in \{0, \ldots, N\}$ and assuming local thermal equilibrium between solid and fluid, energy conservation takes the form
\begin{equation}
    \label{eq:energy_conservation}
    \dd{\specificVolume{}{i}\effective{\density{}{i}\internalEnergy{}{i}}}{t} + \nabla \cdot \left[\specificVolume{}{i} \left( \enthalpyFlux{}{i} + \heatFlux{}{i} \right)\right] - \sum_{j\in\higherSet[i]} \projectToInterface{i}{j} \specificVolume{}{j}\left(\interfaceEnthalpyFlux{}{j} + \interfaceHeatFlux{}{j}\right)= \energySource{}{i}.
\end{equation}
Here, \internalEnergy{}{} denotes internal energy and curly brackets indicate the porosity-weighted sum of a scalar quantity \scalar{}{} arising from the thermal equilibrium assumption:
\begin{equation}
    \label{eq:effective_energy}
    \effective{\scalar{}{}} = \porosity{}{} \scalar{f}{} + (1-\porosity{}{}) \scalar{s}{}.
\end{equation}
Furthermore, \enthalpyFlux{}{i} and \heatFlux{}{i} are enthalpy and heat fluxes within \subdomain{}{i}, respectively, and \interfaceEnthalpyFlux{}{j} and \interfaceHeatFlux{}{j} are their interface counterparts, while \energySource{}{i} is a source or sink of energy.

Ignoring inertial terms, the momentum conservation equation in the matrix subdomain ($d_i=N$) is
\begin{equation}
    \label{eq:momentum_conservation}
    -\nabla \cdot \stress{}{i}
    = \bodyForce{}{i},
\end{equation}
with \stress{}{i} being the total stress tensor and \bodyForce{}{i} body forces.

\subsection{Contact kinematics}
In this subsection, we consider a matrix-fracture pair \subdomain{}{h} and \subdomain{}{l} of dimensions $d_h=N$ and $d_l=N-1$ and denote the two interfaces on either side of \subdomain{}{l} as \interface{}{\plusSide} and \interface{}{\minusSide}, see middle panel of Figure \ref{fig:geometry}.
We define the fracture normal vector \normal{}{l} to coincide with \normal{}{h} on the \plusSide-side and introduce the fracture contact traction, \traction{}{l}, defined according to the direction of \normal{}{l}.
Denoting a generic vector defined on the fracture by \vect, its normal \normalComponent and tangential \tangentialComponent components on \subdomain{}{l} are
\begin{equation}
    \label{eq:fracture_decomposition}
    \iota_{\normalComponent} = \vect \cdot \normal{}{l}, \quad \vect[\tangentialComponent] = \vect - \iota_{\normalComponent} \normal{}{l}.
\end{equation}
We also introduce the jump in interface displacements across \subdomain{}{l}:
\begin{equation}
    \label{eq:displacement_jump}
    \jump{\displacement{}{}} = \projectFromInterface{l}{\minusSide}\displacement{}{\minusSide} - \projectFromInterface{l}{\plusSide}\displacement{}{\plusSide}.
\end{equation}

We achieve balance between the traction on the two fracture surfaces by enforcing each of them to equal the total fracture traction:
\begin{align}
    \label{eq:traction_balance}
    \projectToInterface{h}{\plusSide} \stress{}{h} \cdot \normal{}{h}  & =
    \projectToInterface{l}{\plusSide} \left(\traction{}{l}  - \pressure{}{l}\identity \cdot \normal{}{l} \right),       \\
    -\projectToInterface{h}{\minusSide} \stress{}{h} \cdot \normal{}{h} & =
    \projectToInterface{l}{\minusSide} \left(\traction{}{l}  - \pressure{}{l}\identity\cdot \normal{}{l} \right).
\end{align}

Proceeding to the relations governing the fracture contact kinematics and suppressing the $l$ subscript for readability, non-penetration for the fracture surfaces reads
\begin{align}
    \label{eq:normal_fracture_kinematics}
     \begin{array}{ r l}
    \jump{\displacement{}{}}_{\normalComponent} - \gap{}{} &\geq 0, \\
    \normalTraction  &\leq 0, \\
    \normalTraction \left( \jump{\displacement{}{}}_{\normalComponent} - \gap{}{} \right) & = 0,
    \end{array}
\end{align}
where \gap{}{} is the gap between the fracture surfaces when in mechanical contact.
The second inequality reflects that compressive normal contact traction corresponds to negative \normalTraction by the definition of \traction{}{}.
Denoting the friction bound by \frictionBound and the increment in tangential displacement by \incrementTangentialDisplacement, the friction model is
\begin{align}
    \label{eq:tangential_fracture_kinematics}
     \begin{array}{ r l l l}
                \norm{\traction{}{\tangentialComponent}} &\leq \frictionBound , & & \\
                 \norm{\traction{}{\tangentialComponent}} &< \frictionBound & \rightarrow & \incrementTangentialDisplacement = 0 ,\\
                 \norm{\traction{}{\tangentialComponent}} &=\frictionBound & \rightarrow  & \exists \,\zeta \in \mathbb{R^+}:   \incrementTangentialDisplacement = \zeta \traction{}{\tangentialComponent}.
  \end{array}
\end{align}
The three relations state that: (i) tangential stresses are bounded, (ii) tangential deformation occurs only if the bound is reached, and (iii) tangential stresses and deformation increments are parallel.

\subsection{Constitutive relations}\label{sec:constitutive:relations}
To complement the equations presented in the previous subsections, we select constitutive laws based on Coussy \cite{coussy2004poromechanics}.
Similar models for thermo-poromechanics are found in e.g., \cite{nikolaevskij1990mechanics,garipov2019thm}.

The volumetric fluid flux is modelled using Darcy's law,
\begin{equation}
    \label{eq:darcy}
    \fluidFlux{}{}= - \frac{\permeability{}{}}{\viscosity{}{}} \left(\nabla \pressure{}{}-\density{f}{}\gravityVector\right),
\end{equation}
where \permeability{}{} is the permeability tensor, \viscosity{}{} is fluid viscosity, and \gravityVector the gravitational acceleration vector.
While included in this section for completeness, gravity effects are neglected (equivalent to $\gravityVector=0$) in the simulation examples shown in Sections \ref{sec:testing} and \ref{sec:application:examples}.
We assume \permeability{}{} to be constant in the matrix, whereas the fracture permeability is given by the cubic law:
\begin{equation}
    \permeability{}{i} = \frac{\aperture{2}{i}}{12}, \qquad d_i = N-1. \label{eq:fracture_permeability}
\end{equation}
Note that \aperture{}{} changes depending on \jump{\displacement{}{}}\ as detailed in \eqref{eq:aperture}.
Intersection permeability is computed as the average of the permeability in the intersecting fractures.

Solid density is assumed constant, whereas fluid density is given by
\begin{equation}
    \label{eq:density}
    \density{f}{} = \density{0}{} \exp \left( \compressibility{}{} \left( \pressure{}{} - \pressure{0}{} \right)
    - \thermalExpansion{f}{} \left( \temperature{}{} - \temperature{0}{} \right) \right),
\end{equation}
with \compressibility{}{} denoting compressibility, \thermalExpansion{f}{} fluid thermal expansion coefficient and the superscript $0$ a reference state.
We shall use both a constant viscosity model and a temperature-dependent one, given by
\begin{equation}
    \label{eq:viscosity}
    \viscosity{}{} = \viscosity{A}{} \exp \left( \frac{ \viscosity{B}{}}{\temperature{}{} - \viscosity{C}{}}  \right),
\end{equation}
where \viscosity{A}{}, \viscosity{B}{} and \viscosity{C}{} are fluid constants \cite{viswanath1989data}.

Denoting specific heat capacity by \heatCapacity{}{} and assuming a simple fluid description, fluid enthalpy is given by
\begin{equation}
    \label{eq:enthalpy}
    \enthalpy{f}{} = \heatCapacity{f}{} \left( \temperature{}{} - \temperature{0}{} \right),
\end{equation}
and specific internal energies are computed as
\begin{align}
    \label{eq:internal_energy}
    \internalEnergy{f}{} &= \enthalpy{f}{} - \frac{\pressure{}{}}{\density{f}{}}, \\
    \internalEnergy{s}{} &= \heatCapacity{s}{} \left( \temperature{}{} - \temperature{0}{} \right).
\end{align}

Using the effective thermal conductivity  \effective{\thermalConductivity{}{}},  Fourier's law for the diffusive heat flux through fluid and solid phase reads
\begin{equation}
    \label{eq:fourier}
    \heatFlux{}{} = - \effective{\thermalConductivity{}{}} \nabla \temperature{}{},
\end{equation}
while the advective heat flux is given by
\begin{equation}
    \label{eq:advection}
    \enthalpyFlux{}{} =  \enthalpy{f}{} \density{f}{} \fluidVelocity{}{}.
\end{equation}

The volumetric interface flux is proportional to the pressure jump across \interface{}{j} via a Darcy-type law \cite{Martin2005}:
\begin{equation}
    \label{eq:interfaceFluxes}
    \interfaceFluidFlux{}{j} = - \frac{\permeability{}{j}}{\viscosity{}{j}} \left[ \frac{2}{\projectToInterface{l}{j} \aperture{}{l}} \left(\projectToInterface{l}{j} \pressure{}{l} - \projectToInterface{h}{j} \pressure{}{h}\right) - \gravityVector\density{}{j}\right],
\end{equation}
where \permeability{}{j} is the interface permeability.
Similarly, the diffusive interface heat flux is
\begin{equation}
    \label{eq:interfaceDiffusiveHeatFlux}
    \interfaceHeatFlux{}{j} = - \thermalConductivity{}{j} \frac{2}{\projectToInterface{l}{j} \aperture{}{l}} \left(\projectToInterface{l}{j} \temperature{}{l} - \projectToInterface{h}{j} \temperature{}{h}\right),
\end{equation}
while its advective counterpart is
\begin{equation}
    \label{eq:interfaceEnthalpyFlux}
    \interfaceEnthalpyFlux{}{j} =  \enthalpy{f}{j}\density{f}{j}\interfaceFluidFlux{}{j}.
\end{equation}
The values for \permeability{}{j} and \thermalConductivity{}{j} are both inherited from the lower-dimensional neighbour subdomain.
For an advected quantity \scalar{}{j} representing \density{}{j}, \viscosity{}{j} and \enthalpy{}{j}, we use an inter-dimensional upwinding based on \interfaceFluidFlux{}{j}:
\begin{gather}
\begin{aligned}
    \label{eq:inter_dimensional_upwind}
    \scalar{}{j} = \left\{ \begin{array}{ l l }
      \projectToInterface{h}{j} \scalar{}{h} & \text{ if } \interfaceFluidFlux{}{j} > 0  \\
      \projectToInterface{l}{j} \scalar{}{l} & \text{ if } \interfaceFluidFlux{}{j} \leq 0 .
  \end{array} \right.
\end{aligned}
\end{gather}

The total thermo-poromechanical stress tensor is given by an extended Hooke's law that also accounts for fluid and thermal contributions
\begin{equation}
    \label{eq:total_stress}
   \stress{}{} = \stress{0}{} + \shearModulus(\nabla \displacement{}{}+\nabla \displacement{T}{}{})+(\bulkModulus{}-\frac{2\shearModulus}{3}) \trace{\nabla \displacement{}{}} \identity
 - \thermalExpansion{s}{} \bulkModulus{}{} \left(\temperature{}{}-\temperature{0}{}\right) \identity - \biotAlpha \left(\pressure{}{}-\pressure{0}{}\right) \identity.
\end{equation}
Here, \shearModulus is the shear modulus, \bulkModulus{}{} is the matrix bulk modulus, \biotAlpha is the Biot coefficient, \thermalExpansion{s}{} is the drained thermal expansion coefficient (which equals the solid thermal expansion coefficient) and $\trace{\cdot}$ denotes the trace of a matrix. The assembly is illustrated in Snippet \ref{code:stress}.
Gravitational forces are included in the momentum balance by setting
\begin{equation}
    \label{eq:body_forces}
    \bodyForce{}{} = \effective{\density{}{i}}\gravityVector.
\end{equation}
Matrix porosity depends on pressure, displacement and temperature according to
\begin{equation}
    \label{eq:porosity}
    \porosity{}{} = \porosity{0}{}
    + \biotAlpha \nabla \cdot \displacement{}{}
    + \frac{\left( \biotAlpha-\porosity{0}{} \right) \left( 1 - \biotAlpha \right)}{\bulkModulus{}} \left( \pressure{}{} - \pressure{0}{} \right)
    - \thermalExpansion{\porosity{}{}}{} \left( \temperature{}{} - \temperature{0}{} \right),
\end{equation}
with $\thermalExpansion{\porosity{}{}}{}:=(\biotAlpha-\porosity{0}{})\thermalExpansion{s}{}$ denoting the porosity related thermal expansion coefficient.
We assume unitary fracture and intersection porosity.

\begin{snippet}[ht]
\caption{The thermo-poromechanical stress is assembled by collecting three terms. This allows for seamless code reuse, since the \texttt{stress} method is used with the same signature in the purely mechanical model and the poromechanical to collect the first one and two terms, respectively.
Below, we show the definition of the pressure term, which employs the \texttt{BiotAd} class for discretization of poromechanical terms, as discussed in Section \ref{sec:numerical}.\label{code:stress}}
\fbox{\includegraphics[width=0.975\textwidth]{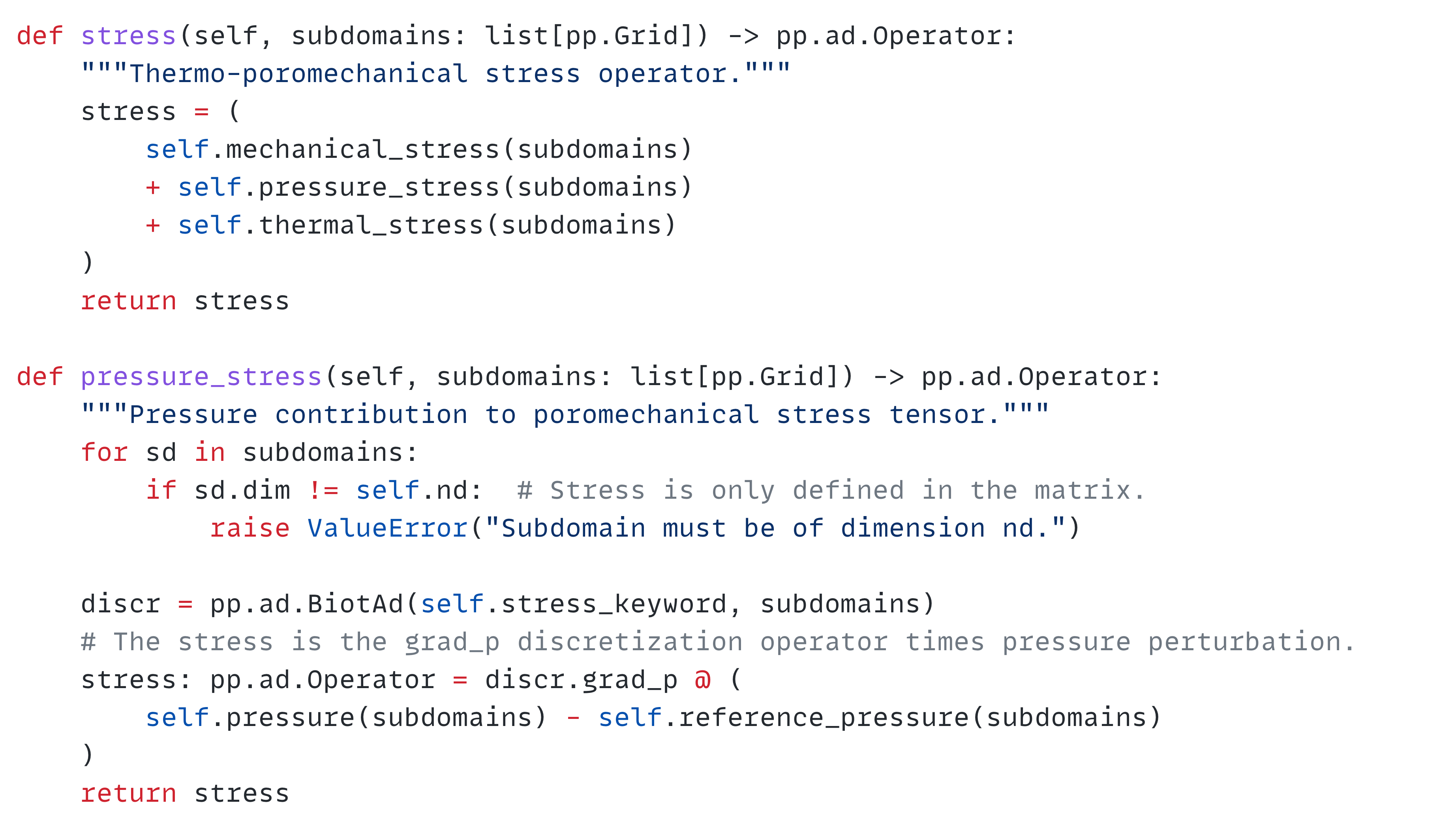}}
\end{snippet}

Turning now to fracture deformation, we assume the friction bound to relate to normal traction through a Coulomb type friction law with a constant friction coefficient \frictionCoefficient:
\begin{equation}
    \label{eq:coulomb_friction}
    \frictionBound = -\frictionCoefficient \normalTraction.
\end{equation}
Fracture roughness effects  are incorporated through the gap function:
\begin{equation}
    \label{eq:fracture_gap}
  \gap{}{} = \gap{0}{} + \tan{\dilationAngle}\norm{\jump{\displacement{}{}}_{\tangentialComponent}} + \frac{\maximumFractureClosure\normalTraction}{\maximumFractureClosure\normalStiffness-\normalTraction}.
\end{equation}
Here, the first term is the residual gap corresponding to an unstressed and undeformed fracture. The second term represents shear dilation, with \dilationAngle denoting the dilation angle.
The third term accounts for elastic normal deformation according to \cite{barton1985strength}, where \normalStiffness is the normal stiffness per area and \maximumFractureClosure is the maximum elastic normal closure of the fracture.
These changes to \gap{}{} impact normal deformation according to Eqs.\ \eqref{eq:normal_fracture_kinematics}, which in turn translate into aperture changes by setting
\begin{equation}
\label{eq:aperture}
\aperture{}{}=\aperture{0}{}+\jump{\displacement{}{}}_{\normalComponent}
\end{equation}
in the fractures.
Here, \aperture{0}{} denotes a residual hydraulic aperture.
In intersection subdomains (i.e., $d_i<N-1$), we compute the aperture as the mean among apertures of neighbouring higher-dimensional subdomains, \aperture{}{h},
\begin{equation}
    \label{eq:intersection_aperture}
    \aperture{}{i} = \frac{1}{|\higherSet[i]|}\sum_{j\in \higherSet[i]}\projectFromInterface{i}{j}\projectToInterface{h}{j}\aperture{}{h}.
\end{equation}
Finally, by interpreting aperture as the side length of the reduced dimension, we obtain the specific volume as
\begin{equation}
    \label{eq:specific_volume}
    \specificVolume{}{i} = \aperture{N-d_i}{i},
\end{equation}
which reduces to unity in the matrix subdomain.

\subsection{Initial and boundary conditions}\label{sec:initial_and_boundary_conditions}

To close the system of equations, we provide initial values for all primary variables and boundary conditions.
For primary variables, we use \pressure{}{i}, \temperature{}{i}, \displacement{}{i} ($d_i=N$) and \traction{}{i} ($d_i=N-1$) in the subdomains and \interfaceFluidFlux{}{j}, \interfaceEnthalpyFlux{}{j}, \interfaceHeatFlux{}{j} and \displacement{}{j} ($d_j=N-1$) on interfaces.

Boundary conditions are set both on internal and external boundaries of each subdomain. On internal boundaries $\partial_j\Omega_i$, we require continuity of normal mass fluxes, normal energy fluxes, and displacement (for $d_i = N$):
\begin{align}
    \specificVolume{}{i} \density{f}{i} \fluidVelocity{}{i} \cdot \normal{}{i} &= \projectFromInterface{i}{j} \specificVolume{}{j} \density{f}{j}\interfaceFluidFlux{}{j}, \label{eq:internal_bc_mass}\\
    \specificVolume{}{i} \left(\enthalpyFlux{}{i} + \heatFlux{}{i}\right) \cdot \normal{}{i} &= \projectFromInterface{i}{j}
    \specificVolume{}{j} \left(
    \interfaceEnthalpyFlux{}{j} + \interfaceHeatFlux{}{j}\right), \label{eq:internal_bc_energy}\\
    \displacement{}{i} &= \projectFromInterface{i}{j}\displacement{}{j}. \label{eq:internal_bc_momentum}
\end{align}
Here, \normal{}{i} is the normal vector on $\partial_j\Omega_i$ pointing from the higher-dimensional to the lower-dimensional subdomain.

On external boundaries, we allow for two types of boundary conditions, namely Neumann and Dirichlet. For boundary conditions of the Neumann type, we prescribe values for mass flux $\specificVolume{}{i} \density{f}{i} \fluidVelocity{}{i} \cdot \normal{}{i}$, energy flux $\specificVolume{}{i} \left(\heatFlux{}{i} + \enthalpyFlux{}{i} \right) \cdot \normal{}{i}$, and total traction $\stress{}{i} \cdot \normal{}{i}$ ($d_i = N$). On immersed fracture tips, we require both mass and energy fluxes to equate to zero.
Note that \normal{}{i} here denotes the outward pointing normal on the exterior boundary.
Finally, for boundary conditions of the Dirichlet type, we prescribe values for pressure \pressure{}{i}, temperature \temperature{}{i}, and displacement \displacement{}{i} ($d_i = N$).

\section{Numerical solution approach \label{sec:numerical}}
The simulation framework \texttt{PorePy} aims to support the mathematical modelling framework described in Section \ref{sec:mathematical:models}.
Practical usage is mainly based on the \texttt{Models} introduced in Section \ref{sec:design:principles}, which are composed of classes defining model equations, geometry, variables, solution strategy, etc.
This modularisation facilitates a high degree of flexibility and code reuse.
Moreover, inspired by the single responsibility principle, each method of these classes performs a limited action.
Combined, this allows for adjustments on multiple levels with minimal effort and intrusiveness.
% - from adding a new balance equation to changing
Below, we detail some numerical aspects of \texttt{PorePy}, including available discretization schemes and treatment of non-linearities by automatic differentiation.

The computational grids are constructed to conform to immersed lower-dimensional geometric objects (fractures and intersections) in the sense that each lower-dimensional subdomain coincides geometrically with a set of faces on the surrounding higher-dimensional grid.
Both logically Cartesian and simplex grids are supported, with the latter generated by \texttt{Gmsh} \cite{geuzaine2009gmsh}.
The data structure for the mixed-dimensional grid contains separate grids for individual subdomains \subdomain{}{i}, mortar grids on the interfaces \interface{}{j}, and discrete versions of the projection operators between subdomains and interfaces.
This structure is also exploited when results are exported for visualization, for which we use \texttt{ParaView} (version 5.11.0 herein).

The choice of discretization methods is motivated by the modelling principles followed in Section \ref{sec:mathematical:models}:
Local conservation is enforced for mass, energy and momentum by applying finite volume methods, specifically, we consider lowest order cell-centred methods.
The discretization of constitutive laws follows their mathematical structure.
Advective terms are discretized by standard upwinding \cite{aziz1979petroleum}.
For diffusive terms, such as those related to Darcy's law in the mass conservation equation and Fourier's law in the energy conservation equation, we use either the standard two-point flux approximation or its more accurate multi-point extension \cite{aavatsmark2002introduction}.
For the stress, we rely solely on the multi-point approximation, which also provides a discretization of the thermo-poromechanical coupling terms, see \cite{nordbotten2016stable,nordbotten2021mpxa}.
The equations on different subdomains are discretized separately, following \cite{nordbotten2019unified, PorePy}, and neighbouring subdomains are coupled via discrete interface variables.

While most of the non-linearities in the governing equations are common for standard reservoir simulation, the contact mechanics problem needs special treatment.
Here we follow techniques from computational contact mechanics:
The inequalities \eqref{eq:normal_fracture_kinematics} and \eqref{eq:tangential_fracture_kinematics} can be reformulated into equivalent equalities as detailed in e.g., \cite{Huber2005,Stefansson2021thm}.
The reformulation allows us to discretize these equations using AD, resulting in a semi-smooth Newton method
\cite{hueber2005primal,berge2020poroelastic_contact}.

Practical usage of \texttt{PorePy} entails experimentation with modelling choices and approaches to simulations for complex and non-linear models in mixed-dimensional geometries.
Accordingly, the code is designed to support rapid prototyping and flexibility in modelling approaches, as was indicated by the snippets in Section \ref{sec:mathematical:models} and further demonstrated in Section \ref{sec:application:examples}.
This is reflected in the data structures underlying the representation and linearization of governing equations:
Provided dimensional compatibility, variables, discretizations and constitutive relations can be defined on arbitrary sets of subdomains and interfaces, thereby enabling tailored governing equations for specific subdomains.
Governing equations are considered non-linear by default, with residual evaluation and linearization implemented using automatic differentiation.
This is again implemented in two steps: Equations are represented as what we term AD operators.
These are symbolic representations of mathematical expressions, which are stored as computational graphs, borrowing popular techniques in e.g., machine learning libraries \cite{paszke2017automatic,abadi2017computational}.
AD operators can be combined by arithmetic operations to form compound AD operators and thereby allow for gradually composing complex expressions.
The translation of the graph into numerical values for the residual and the Jacobian matrix is implemented using forward-mode automatic differentiation \cite{naumann2011art}.

For differentiation in time, we use the backward Euler scheme.
Most of the state-dependent parameters that enter constitutive laws are discretized fully implicitly, that is, their derivatives are included in the Jacobian matrix. The exceptions are permeability, thermal conductivity, and the upwind directions, which enter discretization schemes (multi-point approximations and upwinding) as expressions that are not readily differentiated. These dependencies are lagged one Newton iteration. The high degree of implicitness is made possible by the AD framework which removes the need for manual differentiation of complex expressions.

The following two sections present \texttt{PorePy} examples of testing and multiphysics simulation.
The complete collection of source code, run scripts and simulation results are available as a \texttt{Docker} image \cite{stefansson2023source_code}.

\section{Testing \label{sec:testing}}

Testing must constitute a fundamental part of the development process of scientific software \cite{khan2014importance}. A solid testing framework: (1) reduces the chances of errors by finding problems at an early stage, (2) enhances the quality and reliability of the software and (3) gives the developers a solid ground to incorporate new functionality.

The large number of components present in a multiphysics software demands a systematic and well-structured testing strategy.
According to \cite{burnstein2006practical}, tests can be classified into four levels: (1) unit tests, (2) integration tests, (3) system tests and (4) acceptance tests.
These levels usually form a bottom-up hierarchical structure, with unit tests at the bottom and acceptance tests at the top.

Unit tests check individual components of the software, ideally in isolation.
Integration tests check if groups of individual components interact as expected.
System tests are designed to test if the software as a whole works properly.
Finally, acceptance tests are meant to check whether the software meets the requirements set by an end user.
Since acceptance tests are mostly relevant for commercial software, we will not consider this type of test.

\begin{figure}
    \centering
    \includegraphics[width=\textwidth]{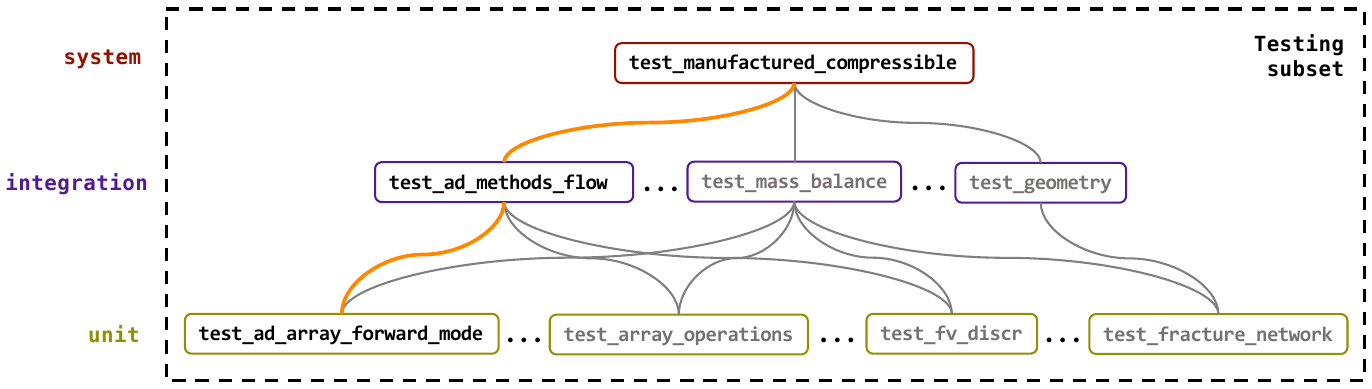}
    \caption{Tests for the compressible single-phase flow model in fractured porous media. There exists a hierarchical structure divided into three levels, i.e., system, integration and unit tests. Note that a small subset of all tests are included in the diagram. Snippets for the tests connected by the orange curves are shown in Sections~\ref{sec:unit}, \ref{sec:integration} and \ref{sec:system}. \label{fig:tests_general}}
\end{figure}

Using the compressible single-phase flow model in fractured porous media\footnote{Recall that this model is given by the conservation law \eqref{eq:mass_conservation}, the constitutive relationships \eqref{eq:darcy} to \eqref{eq:density}, the internal boundary condition \eqref{eq:internal_bc_mass}, and external boundary conditions and initial conditions.}, we will devote the rest of this section to providing concrete examples of unit, integration and system tests.
A schematic representation of a testing subset is shown in Figure~\ref{fig:tests_general}.
For demonstrative purposes, we have chosen three tests from the modules connected by the orange curves.

Most programming languages have libraries that offer testing-specific functionality.
\texttt{PorePy} employs a combination of the testing libraries \texttt{unittest} and \texttt{pytest} \cite{pytest71}.
A particularly attractive feature of \texttt{pytest} is the possibility to parameterise tests.
In the following subsections, the reader will hopefully appreciate that this feature can be used to test a large number of cases in a succinct and effective manner.

\subsection{Unit test \label{sec:unit}}

The implementation of a native framework for defining equations makes thorough testing paramount to ensure that the framework is correctly implemented and that it stays compatible with updates to upstream dependencies such as \texttt{numpy} \cite{harris2020array} and \texttt{scipy} \cite{virtanen2020scipy}.
In this test, we check that AD operators are correctly combined via arithmetic operations (see Snippet \ref{code:unit_test}). The test is parameterised in four dimensions: The first and second dimensions are the left and right operands of the arithmetic operation, respectively. The tested operands are scalars, dense and sparse arrays (i.e., vectors and sparse matrices) and an AD expression with a non-trivial residual and Jacobian matrix.
The third dimension contains the binary operations, namely: sum, subtraction, multiplication, division, exponentiation and array multiplication. The last dimension establishes whether the quantities are represented in forward-mode AD or as a computational graph, as discussed in Section \ref{sec:numerical}.

For any combination of operands and operation, the test evaluates the resulting expression and compares the results to hard-coded known values. Thus, the test validates the implementation of forward-mode AD and, when the expression is represented as a computational graph, the parsing of this graph into numerical values. Not all combinations of operands and operations are permitted.
As an example, \texttt{scipy} currently does not support adding scalars to sparse matrices, thus attempts at evaluation will raise errors.
For these cases, the test checks that the expected errors are raised.
This ensures that changes in upstream dependencies, including new functionality, will be flagged.

The full test consists of 192 individual tests that provide a robust foundation for all valid low-level combinations of arithmetic operations. Most importantly, the inclusion of a non-trivial AD object among the operands means that, by induction, the test gives confidence to the evaluation of complex expressions used in the definition of multiphysics problems.

\begin{snippet}[!h]
\caption{Example of a unit test that tests the combination of fundamental AD operators using standard arithmetic operations. In this and subsequent snippets, we use ellipses ``$\ldots$'' to indicate non-crucial code that has been omitted for the sake of compactness.}
\fbox{\includegraphics[width=0.975\textwidth]{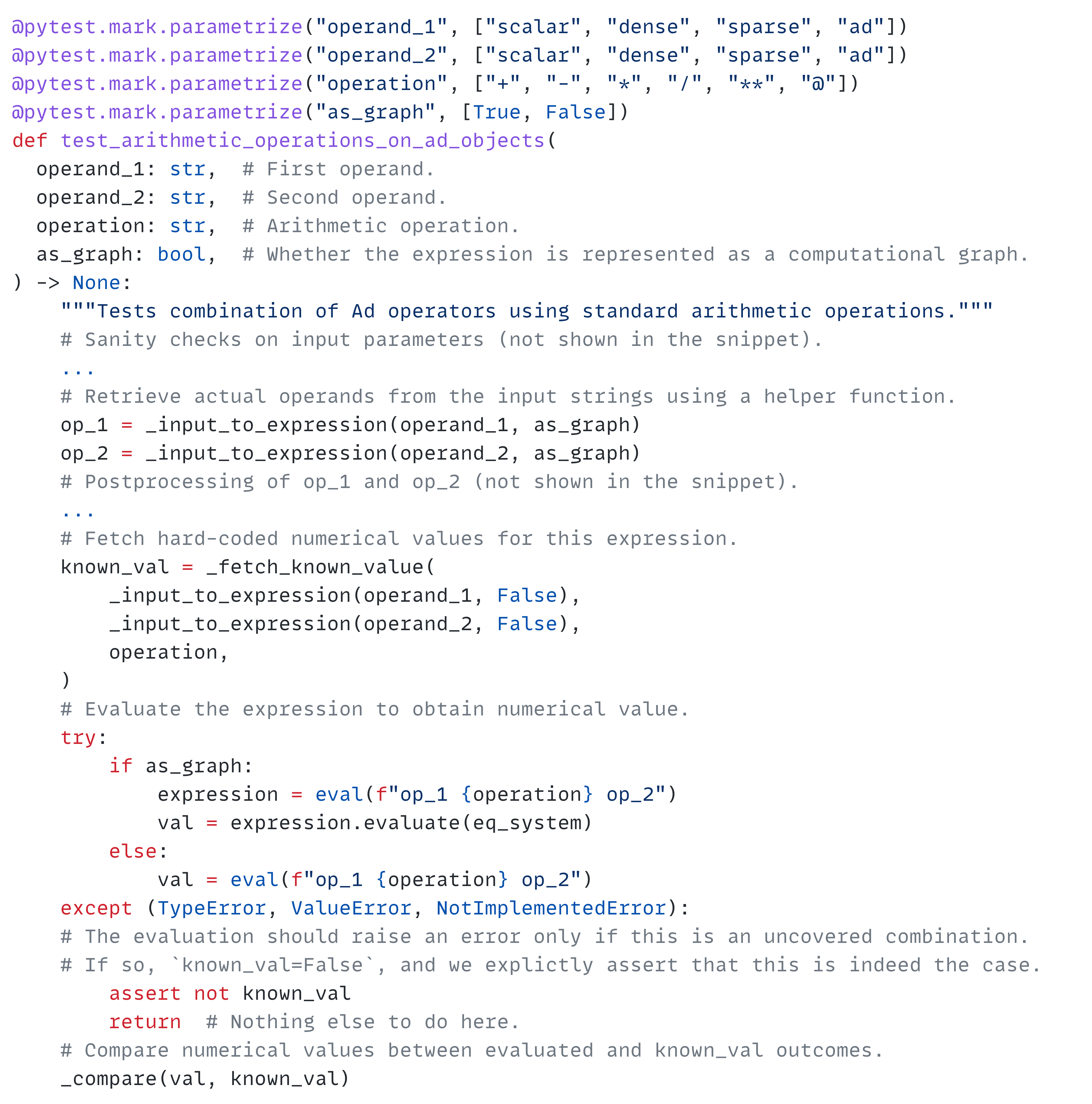}} \label{code:unit_test}
\end{snippet}

\subsection{Integration test\label{sec:integration}}
The modular nature of the \texttt{PorePy} equation definition requires testing not only individual components, but also that the interaction between such components produces expected results.
The integration test shown in Snippet \ref{code:integration_test} compares
the numerical values obtained from methods returning AD operators against known hard-coded values.
The aim of such tests is two-fold: (1) to check that individual methods are correctly implemented and (2) to check that compound expressions obtained from the combination of individual AD operators result in expected values.
Since the process of obtaining a numerical value necessarily requires the resolution of the computational graph and, in some cases, the discretization of relevant quantities, the test is particularly useful to identify parts of the code that are not working.

To minimise the development effort needed to achieve test coverage for new methods, the test requires minimal input related to individual methods: the name of the method, the numerical value expected from evaluation and possible restrictions on which spatial dimensions the method is meaningful.
To cover subdomains and interfaces of different dimensions with a reasonable computational cost, the test is set up on a 2x2 Cartesian grid with two intersecting fractures, which results in subdomains of dimension $0$ to $2$ and interfaces of dimension $0$ and $1$.

The complete test includes $27$ methods.
However, for the sake of compactness, we include only the fluid viscosity and the fluid density as given by Eq.\ \eqref{eq:density}.
The implementation of constant viscosity exemplifies a standalone method, whereas the implementation of the fluid density, which depends on the reference pressure, reference density, compressibility and current pressure state, represents a case where various methods are required to integrate correctly.

\begin{snippet}[!h]
\caption{Example of an integration test that tests the evaluation of AD methods for the compressible single-phase flow model.}
\fbox{\includegraphics[width=0.975\textwidth]{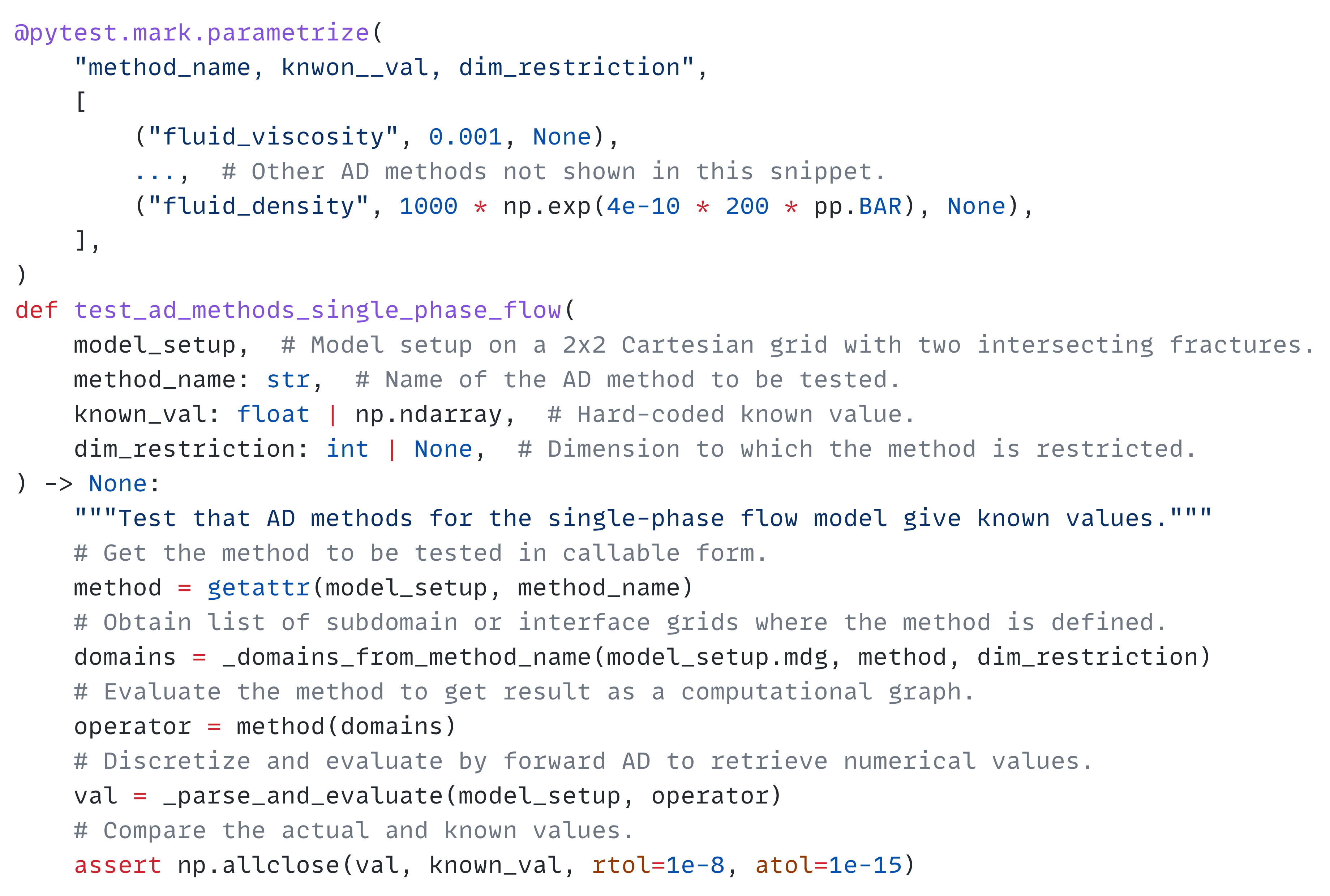}} \label{code:integration_test}
\end{snippet}

\subsection{System test\label{sec:system}}

Ultimately, we would like to know if a model produces the correct results.
Thus, unsurprisingly, system tests are considered the most important type of test in a testing framework \cite{kempf2017system}.
In this context, one of the most robust tests for numerical code verification is the method of manufactured solutions \cite{roy2005review, oberkampf2010verification}.
If available, synthetic solutions represent an invaluable asset for a computational model, as they give confidence that essential parts of the code work as expected.

In a testing context, synthetic solutions are particularly useful when combined with a convergence analysis.
This gives not only a qualitative certainty but also a quantitative validation that errors decrease with increasing spatial and temporal resolution,
thus also verifying implementation of discretization schemes and the like.
More importantly, convergence rates in the asymptotic range are generally invariant to rounding errors and even to minor changes in meshes.
The latter is important for \texttt{PorePy}, since simplex meshes generated with \texttt{Gmsh} may change slightly as this external dependency is updated.

In Snippet \ref{code:functional_test}, we show a test that compares known and actual observed order of converge (OOC) values for the compressible single-phase flow model.
The manufactured solution was obtained generalising the one proposed in \cite{varela2022aposteriori} from the incompressible to the compressible case. The model includes a single fully embedded vertical fracture and gravity effects are neglected. An in-depth explanation of the derivation of the solution is given in \ref{sec:exact_sol}.
Due to the inherent complexity associated to setting up the different cases, performing the convergence analysis and computing the OOC, the test relies on the \texttt{pytest} fixture functions \texttt{desired_ooc} and \texttt{actual_ooc} to collect the known and actual order of convergence rates.

The test is parameterised in three dimensions, namely: variable, grid type and dimension. This variability is what makes this test truly a system test. We check the OOC for primary and secondary variables (matrix pressure, matrix flux, fracture pressure, fracture flux and interface flux) on two types of grids (Cartesian and simplicial) and for two dimensions ($2$ and $3$).

It should be mentioned that spatio-temporal convergence tests are generally resource-intensive, especially in 3D, and we have therefore created independent workflows to run these tests less frequently, i.e., once or twice per week.
We reiterate the importance of the test: It covers all parts of a simulation model, including meshing, discretizing conservation and constitutive laws and boundary and initial conditions, and thereby gives confidence that all these parts are correctly and consistently implemented.

\begin{snippet}[!h]
\caption{Example of a system test that tests observed order of convergence (for primary and secondary variables) for the compressible flow model with a single vertical fracture. These values were used to produce the plot from the right panel of Fig.~\ref{fig:compressible_flow}.\label{code:functional_test}}
\fbox{\includegraphics[width=0.975\textwidth]{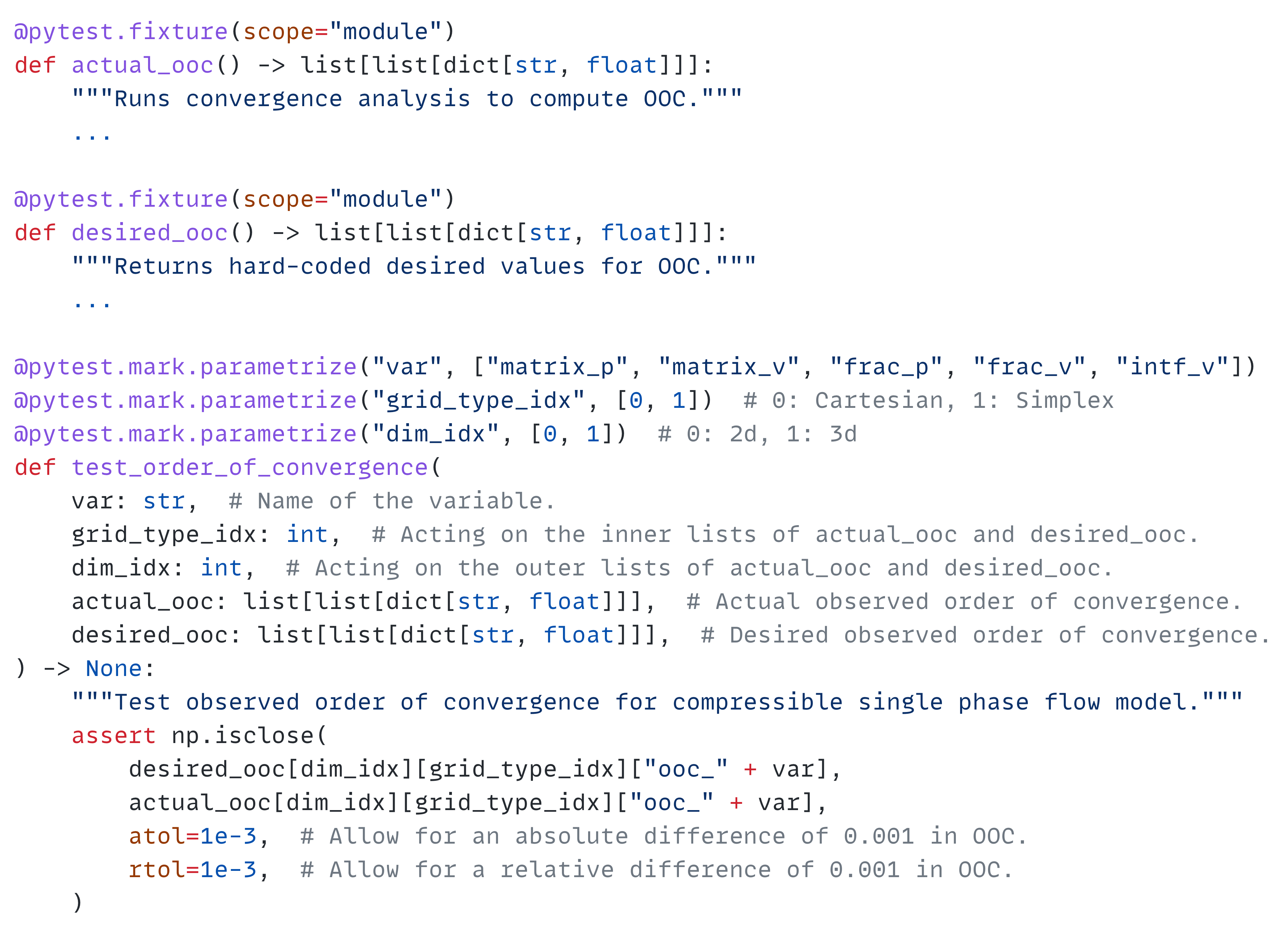}}
\end{snippet}

\section{Application examples}
\label{sec:application:examples}
The two simulations shown in this section illustrate the versatility of \texttt{PorePy} in terms of physical modelling capability.
We stress that the simulations are intended to serve as illustration of \texttt{PorePy} as a modelling tool rather than being interpreted for their physical implications.
For examples of how \texttt{PorePy} is used to study multiphysics processes and develop simulation technology, see e.g. \cite{Stefansson2021thm,Stefansson2023numerical,banshoya2023,stefansson2021propagation,dang2022poro_wing_cracks}.

\subsection{Problem definition\label{sec:application:examples:definition}}
In the first example (Example 1), we solve a mixed-dimensional poromechanical problem with contact mechanics at the fracture interfaces.
This corresponds to Eqs.\ \eqref{eq:mass_conservation}, \eqref{eq:momentum_conservation}-\eqref{eq:tangential_fracture_kinematics}, \eqref{eq:internal_bc_mass} and \eqref{eq:internal_bc_momentum} complemented by the relevant constitutive equations defined in Section \ref{sec:constitutive:relations} (ignoring temperature terms in Eqs.\ \eqref{eq:total_stress} and \eqref{eq:porosity}).
For the second example (Example 2), we extend the model by adding the energy equations \eqref{eq:energy_conservation} and \eqref{eq:internal_bc_energy} and the rest of the constitutive laws.

The domain is a 3D box of dimensions $\SI{100}{\meter}\times\SI{50}{\meter}\times\SI{50}{\meter}$ containing two horizontal fractures, as shown in Figure \ref{fig:applications_pressure_and_temperature}.
The fractures are offset in the vertical direction and partly overlap in the horizontal direction.
Fracture 1 extends to the left domain boundary ($x=\SI{0}{\meter}$) and Fracture 2 extends to the right boundary ($x=\SI{100}{\meter}$), resulting in a geometry which is symmetric about the plane $y=\SI{25}{\meter}$.

We set zero initial displacement and contact traction values and impose Dirichlet displacement boundary conditions on the bottom boundary of the domain and Neumann conditions elsewhere.
On the top, the traction values are \SI{-1e6}{\pascal} and \SI{-2e6}{\pascal} in the $x$ and $z$ direction, respectively.
Combined with zero traction values on the remaining boundaries, these values result in compression and shear displacement.
On the two fracture boundaries, we impose Dirichlet boundary conditions for pressure and temperature.
Initially, the values are $p=\SI{1e5}{\pascal}$ and $T=\SI{400}{\kelvin}$ on both fracture boundaries, thus matching the initial and reference values.
After an equilibration phase of \SI{1.0}{\year},
we change the values on the leftmost (inlet) boundary, increasing to $p=$ \SI{5e5}{\pascal} and reducing to $T=$ \SI{390}{\kelvin}.
On the 3D domain boundaries, we impose homogeneous Neumann boundary conditions for all fluid and energy flux variables.
This produces a flow field from the leftmost fracture, through the central part of the matrix and entering the rightmost fracture before reaching the outlet boundary.

\subsection{Simulation setup}\label{sec:application:example:setup}
As mentioned above, both poromechanics and thermo-poromechanics are among \texttt{PorePy}'s suite of ready-to-run model classes, which contain balance equations, constitutive laws and discretization schemes, default values for material parameters, etc.
Snippet \ref{code:thermo_poro_mechanics} shows how we combine run-script classes implementing the problem specifications described in the previous paragraph.
We illustrate one such class in Snippet \ref{code:geometry}, which shows how we define the mixed-dimensional geometry.
The difference between the two setups indicated in Snippet \ref{code:thermo_poro_mechanics}  illustrate how to adjust both conservation equations and constitutive laws.
The few lines of code needed to override \texttt{PorePy}'s constant viscosity model to the temperature-dependent Eq.\ \eqref{eq:viscosity} is shown in Snippet \ref{code:viscosity_model}.
All material parameter values are listed in Table \ref{tab:parameters}.

\begin{snippet}[!h]
\caption{Definition of a tailored thermo-poromechanical model.
    Each of the collected classes contains implementation of functionality corresponding to its name.
    Differences to poromechanical model are highlighted.
    The class \texttt{PostProcessing} provides functionality for collecting data used to produce the figures in Section \ref{sec:application:examples:results} and is not part of the model setup as such.\label{code:thermo_poro_mechanics}}
\fbox{\includegraphics[width=0.975\textwidth]{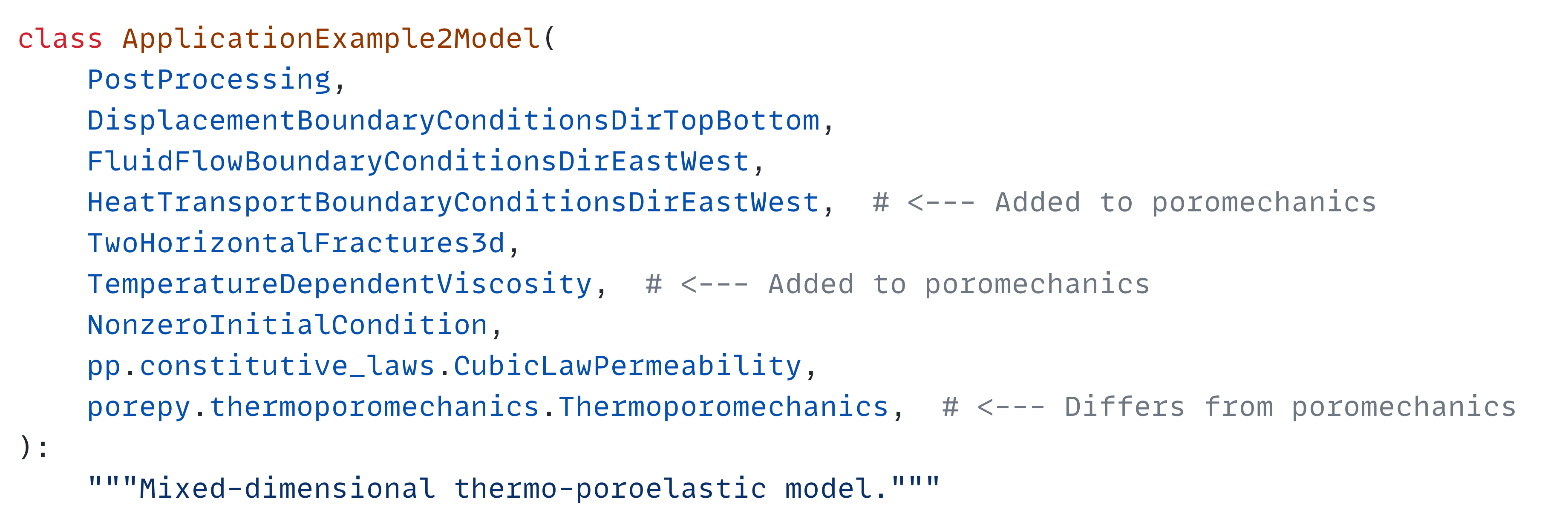}}
\end{snippet}

\begin{table}[ht]
    \centering
    \caption{Material parameters for Section \ref{sec:application:examples}.\label{tab:parameters}}
    \begin{tabular}{|l|l|}
        \hline
        \textbf{Parameter} & \textbf{Value} \\
        \hline
        Biot coefficient, $\alpha$ & \num{8.00e-01} \\
        Matrix permeability, \permeability{}{} & \SI{2.00e-15}{\meter \squared} \\
        Shear modulus, \shearModulus & \SI{1.67e+10}{\pascal} \\
        Bulk modulus, \bulkModulus &  \SI{2.22e+10}{\pascal} \\
        Solid specific heat capacity, \heatCapacity{s}{} & \SI{7.90e+02}{\joule\per \kilogram \per \kelvin} \\
        Solid thermal conductivity, \thermalConductivity{s}{} & \SI{2.50e+00}{\watt\per\meter\per\kelvin} \\
        Solid thermal expansion, \thermalExpansion{s}{} & \SI{1.00e-04}{\per\kelvin} \\
        Fluid specific heat capacity, \heatCapacity{f}{} & \SI{4.18e+03}{\joule\per \kilogram \per \kelvin} \\
        Fluid thermal conductivity, \thermalConductivity{f}{}& \SI{6.00e-1}{\watt\per\meter\per\kelvin} \\
        Fluid thermal expansion, \thermalExpansion{f}{} & \SI{2.10e-04}{\per\kelvin} \\
        Fluid compressibility, \compressibility{}{} & \SI{4e-10}{\per\pascal} \\
        Fluid viscosity (1st simulation), \viscosity{}{} & \SI{1e-3}{\pascal \second} \\
        Viscosity parameter, \viscosity{A}{} & \SI{2.94e-5}{\pascal\second} \\
        Viscosity parameter, \viscosity{B}{} & \SI{5.08e2}{\kelvin} \\
        Viscosity parameter, \viscosity{C}{} & \SI{1.49e2}{\kelvin} \\
        Residual aperture, \aperture{0}{} & \SI{5.00e-4}{\meter} \\
        Solid density, \density{s}{} & \SI{2.70e+03}{\kilogram\per\meter\cubed} \\
        Reference fluid density, \density{0}{} & \SI{1.00e+03}{\kilogram\per\meter\cubed} \\
        Maximum fracture closure, \maximumFractureClosure & \SI{5.00e-4}{\meter} \\
        Fracture normal stiffness, \normalStiffness & \SI{1.00e+09}{\pascal\per \meter} \\
        Residual fracture gap, \gap{0}{} & \SI{5.00e-4}{\meter} \\
        Friction coefficient, \frictionCoefficient & \num{1.00e0} \\
        Dilation angle, \dilationAngle & \SI{5.00e-02}{\radian} \\
        Reference temperature, \temperature{0}{} & \SI{4.00e2}{\kelvin} \\
        Reference pressure, \pressure{0}{} & \SI{1.01e5}{\pascal} \\
        Reference porosity, \porosity{0}{} &\num{5e-2} \\
                \hline
    \end{tabular}
\end{table}

\begin{snippet}[!h]
\caption{Implementation of domain specification, fracture geometry and meshing parameters.
    The adjustments shown herein suffice to produce the geometry shown in Figure \ref{fig:applications_pressure_and_temperature}.
    The call to \texttt{convert_units} ensures consistency if the simulation is run with scaled units to aid numerical performance.
    \label{code:geometry}}
\fbox{\includegraphics[width=0.975\textwidth]{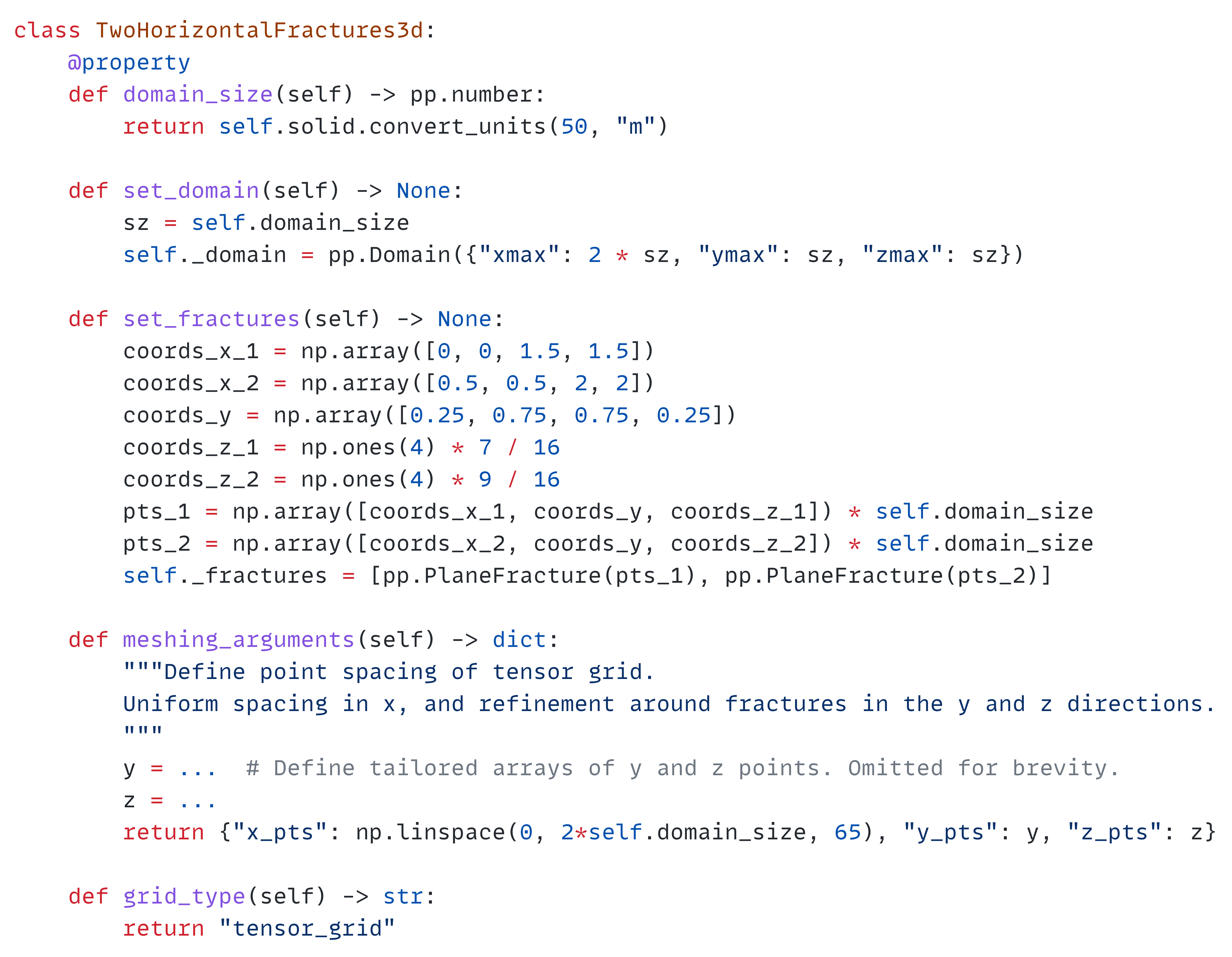}}
\end{snippet}

\begin{snippet}[!h]
\caption{Implementation of a viscosity model in \texttt{PorePy}.\label{code:viscosity_model}}
\fbox{\includegraphics[width=0.975\textwidth]{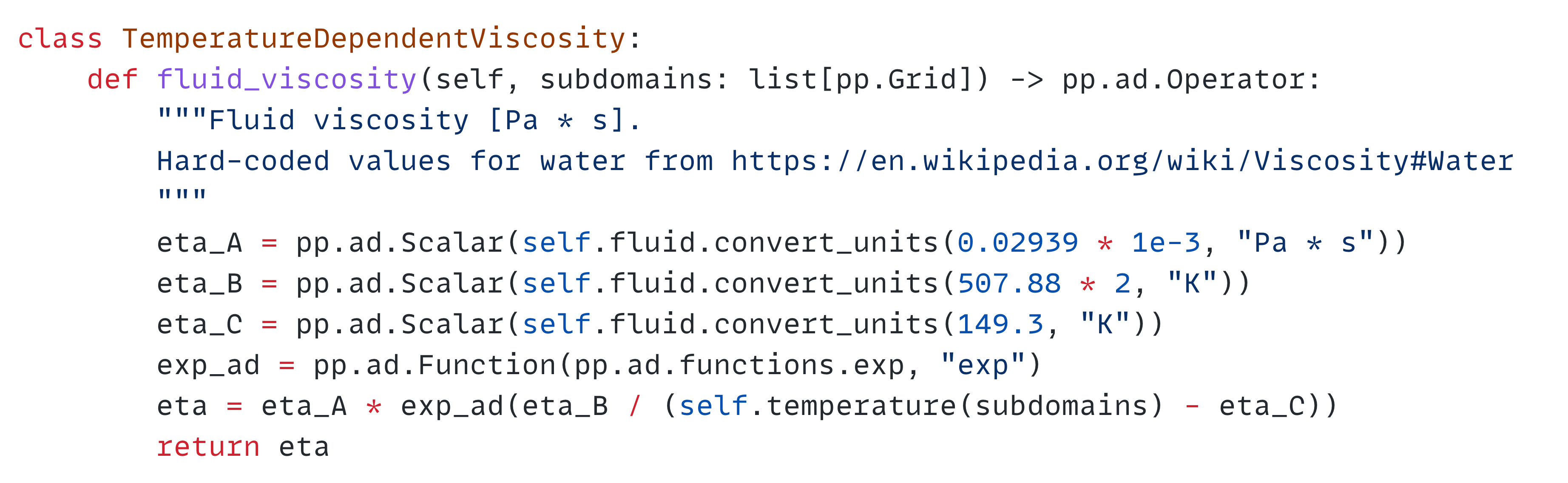}}
\end{snippet}

\subsection{Results}
\label{sec:application:examples:results}
\begin{figure}
    \centering
    \includegraphics[width=0.49\textwidth]{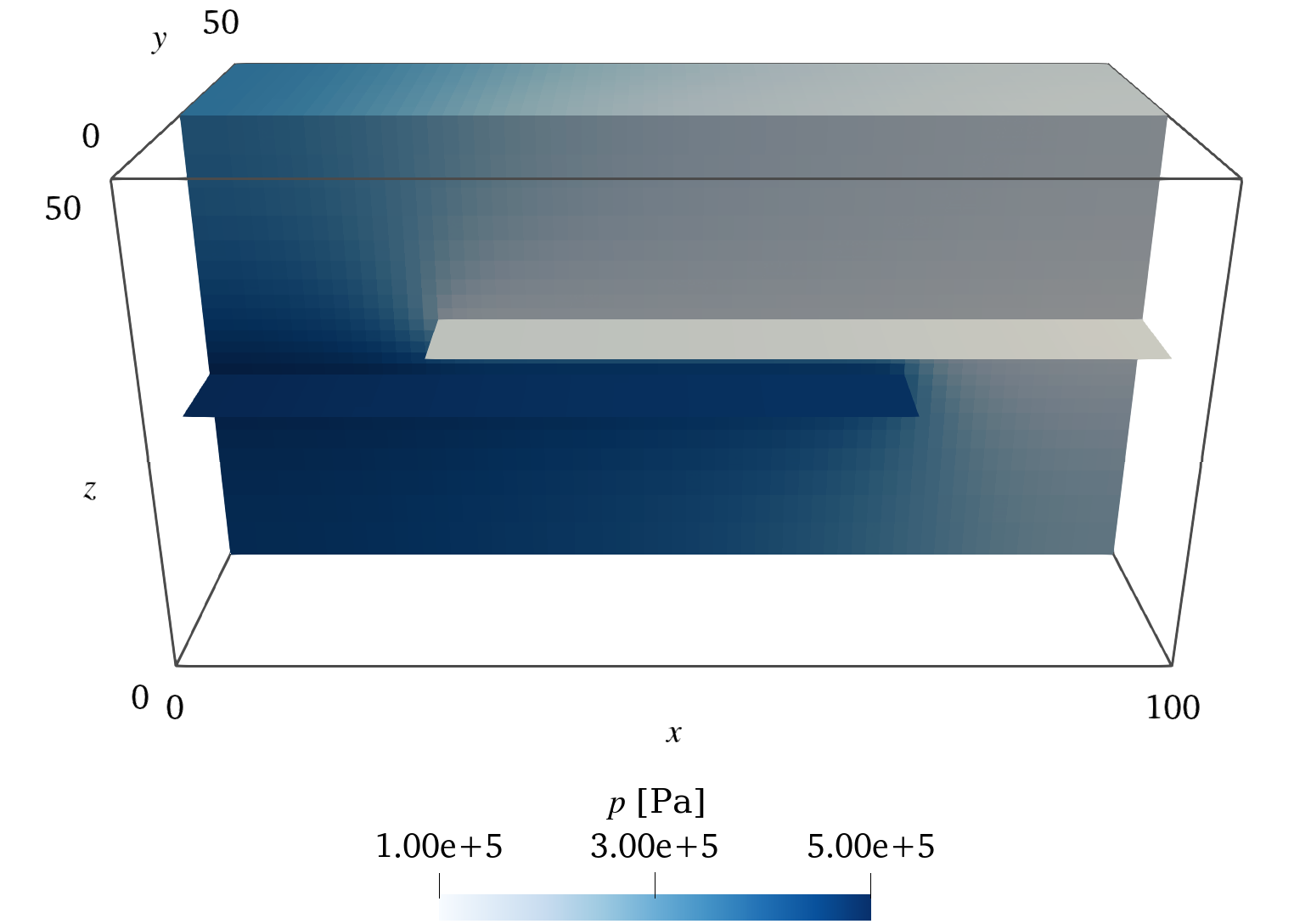}
    \includegraphics[width=0.49\textwidth]{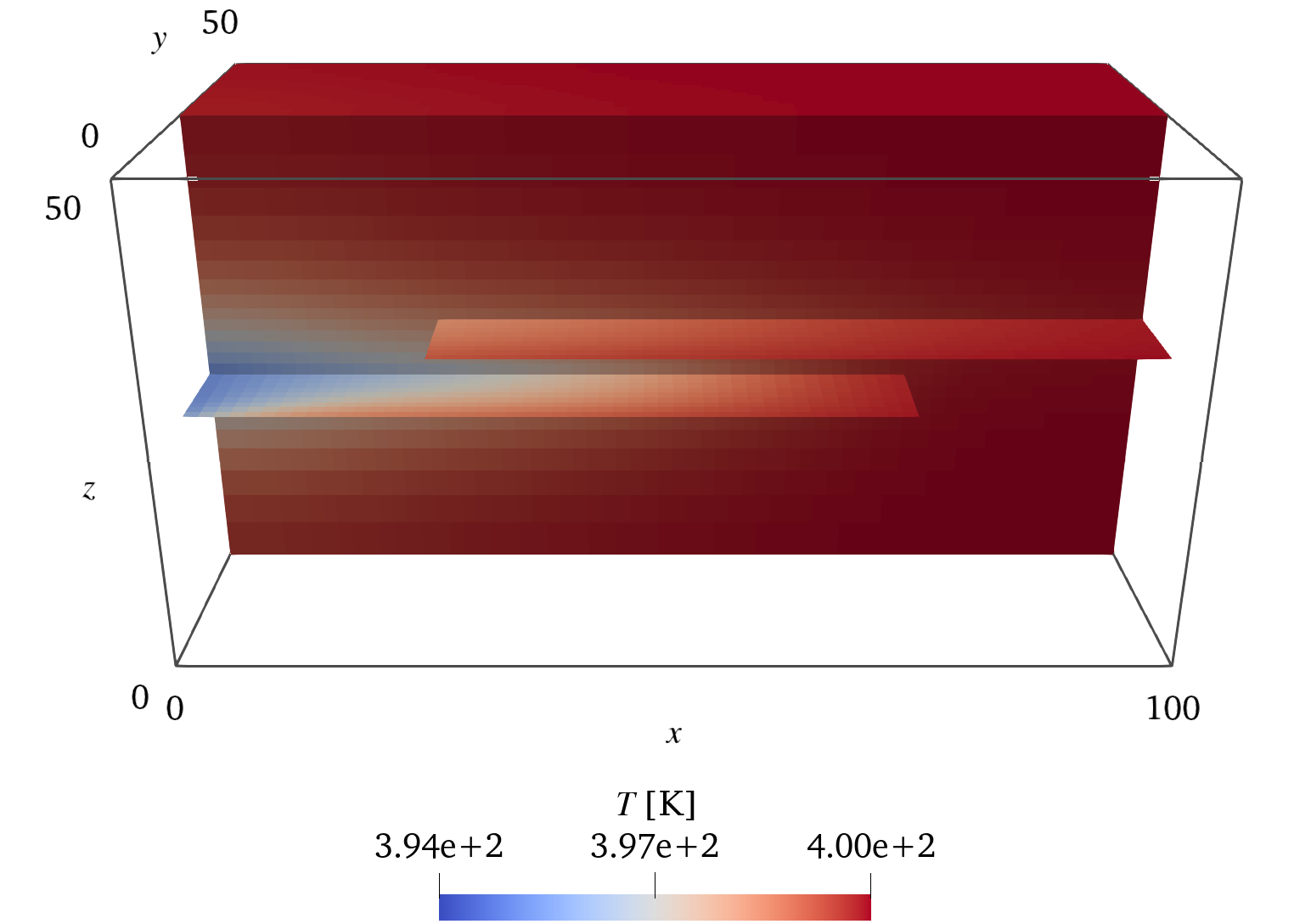}
    \caption{Solutions for the two examples of Section \ref{sec:application:examples}.
    Final pressure distribution for Example 1 (left) and final temperature distribution for Example 2 (right).
    The two lines along the fractures in the leftmost figure indicate the sampling location of the plots shown in Fig. \ref{fig:spatiotemporal_apertures}.}
    \label{fig:applications_pressure_and_temperature}
\end{figure}

Figure \ref{fig:applications_pressure_and_temperature} shows spatial distribution of the pressure and temperature in both fractures and matrix.
The figure underlines the importance of employing separate representations of fractures and matrix to accurately model thermo-poromechanics in fracture media:
The figure shows the fractures to be the main conduit for the pressure pulse and the cooling front.
Moreover, capturing matrix-fracture interaction is also important, as can be seen for both the pressure perturbation and the temperature front.

\begin{figure}
    \centering
    \includegraphics[width=0.98\textwidth]{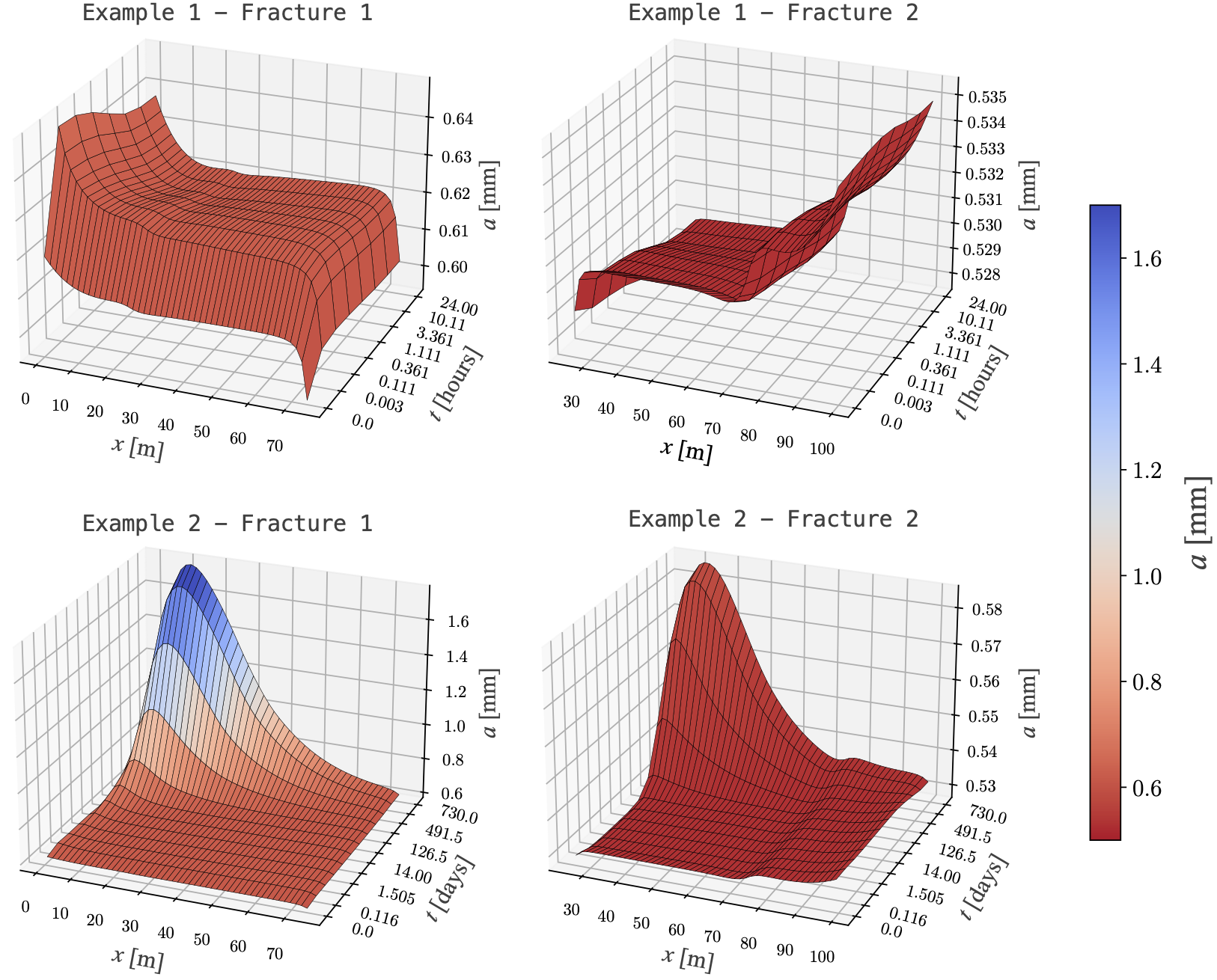}
    \caption{Apertures along central line of the two fractures throughout the first simulation (top) and the second simulation (bottom).
    Note that the time scale for Example 1 and Example 2 are different and they are non-linear of a quasi-logarithmic type.}
    \label{fig:spatiotemporal_apertures}
\end{figure}

Next, consider the temporal evolution of the aperture profiles along the middle of the fractures, shown in Figure \ref{fig:spatiotemporal_apertures}.
The figure contains results of two simulations in two different subdomains along both a spatial and temporal axis, allowing for comparison across four different dimensions.

For the first simulation, the aperture of Fracture 1 can be seen to be largest close to the inlet and decay along the fracture.
This can be attributed to the fluid pressure being highest, thus the effective contact traction smallest, in the vicinity of the inlet.
The effect abates somewhat with time as the pressure signal diffuses and thus illustrates the value of spatiotemporal resolution.
For the aperture in Fracture 2, the most prominent characteristic is the increase at the point where the overlap between Fractures 1 and 2 ends ($x=\SI{75}{\meter}$), illustrating the mechanical coupling between deformation of non-intersecting fractures.

For the second, thermo-poromechanical, simulation, the aperture increase is significantly larger than in the pure poromechanical simulation.
The effect increases with time and is most pronounced close to the inlet (Fracture 1).
It can also be seen in Fracture 2, albeit at a smaller magnitude.
We attribute the difference between the two simulations to thermal contraction due to matrix cooling.

Taken together, the Figure \ref{fig:spatiotemporal_apertures} plots illustrate the type of investigations \texttt{PorePy} facilitates with spatial and temporal analysis of various thermo-poromechanical effects in subdomains of varying dimensions.

\section{Concluding remarks\label{sec:conclusion}}
This paper discusses the design of a multiphysics simulation framework,
presented in the context of thermo-poromechanics in fractured porous media.
Research exploring various coupled processes and constitutive models, and related simulation technology, require flexible software which can be tailored to the governing mathematical model. %A key feature is a code design with a clear correspondence between the governing equations and code.
The simulation toolbox \texttt{PorePy} is structured to closely adhere to the mixed-dimensional governing equations.
The toolbox provides automatic mesh generation and fully coupled discretizations, with spatial derivatives handled by finite volume methods.
\texttt{PorePy} contains a number of  application-relevant multiphysics problems represented by ready-to-run \texttt{Models}.
Thanks to a modular code architecture and automatic differentiation, the \texttt{Models} can be extended  for easy exploration of e.g., constitutive relations and solution strategies with minimal programming.
The reliability of \texttt{PorePy} is ensured by a test suite which monitors the code integrity on a unit, integration and system level.

\section*{Acknowledgements}
This project has received funding from the European Research Council (ERC) under the European Union’s Horizon 2020 research and innovation programme (grant agreement No 101002507).
The work of E. Keilegavlen was financed in part by Norwegian Research Council Grant 308733.

% Bibliographies can be prepared with BibTeX using amsplain.
\bibliographystyle{unsrtnat}
\bibliography{flexible_and_rigorous_arxiv.bib}

\begin{thebibliography}{49}
\providecommand{\natexlab}[1]{#1}
\providecommand{\url}[1]{\texttt{#1}}
\expandafter\ifx\csname urlstyle\endcsname\relax
  \providecommand{\doi}[1]{doi: #1}\else
  \providecommand{\doi}{doi: \begingroup \urlstyle{rm}\Url}\fi

\bibitem[Matth{\"a}i et~al.(2001)Matth{\"a}i, Geiger, and Roberts]{matthai2001csmp}
SK~Matth{\"a}i, S~Geiger, and SG~Roberts.
\newblock Complex systems platform: Csp3d3. 0: user's guide.
\newblock Technical report, ETH Zurich, 2001.

\bibitem[Lie(2019)]{lie2019mrst_book}
Knut-Andreas Lie.
\newblock \emph{An introduction to reservoir simulation using MATLAB/GNU Octave: User guide for the MATLAB Reservoir Simulation Toolbox (MRST)}.
\newblock Cambridge University Press, 2019.

\bibitem[Jung et~al.(2017)Jung, Pau, Finsterle, and Pollyea]{jung2017tough3}
Yoojin Jung, George Shu~Heng Pau, Stefan Finsterle, and Ryan~M Pollyea.
\newblock Tough3: A new efficient version of the tough suite of multiphase flow and transport simulators.
\newblock \emph{Computers \& Geosciences}, 108:\penalty0 2--7, 2017.

\bibitem[Gross and Mazuyer(2021)]{gross2021geosx}
Herve Gross and Antoine Mazuyer.
\newblock Geosx: A multiphysics, multilevel simulator designed for exascale computing.
\newblock In \emph{SPE Reservoir Simulation Conference}, page D011S010R007. OnePetro, 2021.
\newblock \doi{10.2118/203932-MS}.

\bibitem[Podgorney et~al.(2021)Podgorney, Finnila, Simmons, and McLennan]{podgorney2021falcon}
Robert Podgorney, Aleta Finnila, Stuart Simmons, and John McLennan.
\newblock A reference thermal-hydrologic-mechanical native state model of the utah forge enhanced geothermal site.
\newblock \emph{Energies}, 14\penalty0 (16):\penalty0 4758, 2021.

\bibitem[Cacace and Jacquey(2017)]{cacace2017flexible}
Mauro Cacace and Antoine~B Jacquey.
\newblock Flexible parallel implicit modelling of coupled thermal--hydraulic--mechanical processes in fractured rocks.
\newblock \emph{Solid Earth}, 8\penalty0 (5):\penalty0 921--941, 2017.

\bibitem[Koch et~al.(2021)Koch, Gl{\"a}ser, Weishaupt, Ackermann, Beck, Becker, Burbulla, Class, Coltman, Emmert, et~al.]{koch2021dumux}
Timo Koch, Dennis Gl{\"a}ser, Kilian Weishaupt, Sina Ackermann, Martin Beck, Beatrix Becker, Samuel Burbulla, Holger Class, Edward Coltman, Simon Emmert, et~al.
\newblock Dumux 3--an open-source simulator for solving flow and transport problems in porous media with a focus on model coupling.
\newblock \emph{Computers \& Mathematics with Applications}, 81:\penalty0 423--443, 2021.

\bibitem[Bilke et~al.(2022)Bilke, Fischer, Naumov, Lehmann, Wang, Lu, Meng, Rink, Grunwald, Buchwald, Silbermann, Habel, Günther, Mollaali, Meisel, Randow, Einspänner, Shao, Kurgyis, Kolditz, and Garibay]{ogs:6.4.3}
Lars Bilke, Thomas Fischer, Dmitri Naumov, Christoph Lehmann, Wenqing Wang, Renchao Lu, Boyan Meng, Karsten Rink, Norbert Grunwald, Jörg Buchwald, Christian Silbermann, Robert Habel, Linda Günther, Mostafa Mollaali, Tobias Meisel, Jakob Randow, Sophia Einspänner, Haibing Shao, Kata Kurgyis, Olaf Kolditz, and Jaime Garibay.
\newblock Opengeosys, April 2022.

\bibitem[Voskov(2017)]{voskov2017darts}
Denis~V. Voskov.
\newblock Operator-based linearization approach for modeling of multiphase multi-component flow in porous media.
\newblock \emph{Journal of Computational Physics}, 337:\penalty0 275--288, 2017.
\newblock ISSN 0021-9991.
\newblock \doi{https://doi.org/10.1016/j.jcp.2017.02.041}.

\bibitem[B{\v{r}}ezina and Stebel(2016)]{bvrezina2016flow123d}
Jan B{\v{r}}ezina and Jan Stebel.
\newblock Analysis of model error for a continuum-fracture model of porous media flow.
\newblock In \emph{High Performance Computing in Science and Engineering: Second International Conference, HPCSE 2015, Sol{\'a}{\v{n}}, Czech Republic, May 25-28, 2015, Revised Selected Papers 2}, pages 152--160. Springer, 2016.

\bibitem[White et~al.(2018)White, Fu, McClure, Danko, Elsworth, Sonnenthal, Kelkar, and Podgorney]{white2018benchmark_suite}
Mark White, Pengcheng Fu, Mark McClure, George Danko, Derek Elsworth, Eric Sonnenthal, Sharad Kelkar, and Robert Podgorney.
\newblock A suite of benchmark and challenge problems for enhanced geothermal systems.
\newblock \emph{Geomechanics and Geophysics for Geo-Energy and Geo-Resources}, 4:\penalty0 79--117, 2018.

\bibitem[Mindel et~al.(2021)Mindel, Alt-Epping, Landes, Beernink, Birdsell, Bloemendal, Hamm, Lopez, Maragna, Nielsen, Olivella, Perreaux, Saaltink, Saar, {Van den Heuvel}, Vidal, and Driesner]{mindel2021benchmark}
Julian~E. Mindel, Peter Alt-Epping, Antoine Armandine~Les Landes, Stijn Beernink, Daniel~T. Birdsell, Martin Bloemendal, Virginie Hamm, Simon Lopez, Charles Maragna, Carsten~M. Nielsen, Sebastia Olivella, Marc Perreaux, Maarten~W. Saaltink, Martin~O. Saar, Daniela {Van den Heuvel}, Rubén Vidal, and Thomas Driesner.
\newblock Benchmark study of simulators for thermo-hydraulic modelling of low enthalpy geothermal processes.
\newblock \emph{Geothermics}, 96:\penalty0 102130, 2021.
\newblock ISSN 0375-6505.
\newblock \doi{10.1016/j.geothermics.2021.102130}.

\bibitem[Keilegavlen et~al.(2020)Keilegavlen, Berge, Fumagalli, Starnoni, Stefansson, Varela, and Berre]{PorePy}
Eirik Keilegavlen, Runar Berge, Alessio Fumagalli, Michele Starnoni, Ivar Stefansson, Jhabriel Varela, and Inga Berre.
\newblock Porepy: an open-source software for simulation of multiphysics processes in fractured porous media.
\newblock \emph{Computational Geosciences}, 25\penalty0 (1):\penalty0 243--265, 2020.
\newblock ISSN 1420-0597 1573-1499.
\newblock \doi{10.1007/s10596-020-10002-5}.

\bibitem[Martin et~al.(2005)Martin, Jaffré, and Roberts]{Martin2005}
Vincent Martin, Jérôme Jaffré, and Jean~E. Roberts.
\newblock Modeling fractures and barriers as interfaces for flow in porous media.
\newblock \emph{SIAM Journal on Scientific Computing}, 26\penalty0 (5):\penalty0 1667--1691, 2005.
\newblock ISSN 1064-8275 1095-7197.
\newblock \doi{10.1137/s1064827503429363}.

\bibitem[Boon et~al.(2018)Boon, Nordbotten, and Yotov]{boon2018robust}
Wietse~M Boon, Jan~M Nordbotten, and Ivan Yotov.
\newblock Robust discretization of flow in fractured porous media.
\newblock \emph{SIAM Journal on Numerical Analysis}, 56\penalty0 (4):\penalty0 2203--2233, 2018.

\bibitem[Nordbotten et~al.(2019)Nordbotten, Boon, Fumagalli, and Keilegavlen]{nordbotten2019unified}
J.~M. Nordbotten, W.~M. Boon, A.~Fumagalli, and E.~Keilegavlen.
\newblock Unified approach to discretization of flow in fractured porous media.
\newblock \emph{Computational Geosciences}, 23\penalty0 (2):\penalty0 225--237, 2019.
\newblock ISSN 1420-0597 1573-1499.
\newblock \doi{10.1007/s10596-018-9778-9}.

\bibitem[Boon et~al.(2021)Boon, Nordbotten, and Vatne]{boon2021mixed_dimensional}
Wietse~M Boon, Jan~M Nordbotten, and Jon~E Vatne.
\newblock Functional analysis and exterior calculus on mixed-dimensional geometries.
\newblock \emph{Annali di Matematica Pura ed Applicata (1923-)}, 200\penalty0 (2):\penalty0 757--789, 2021.

\bibitem[Boon and Nordbotten(2022)]{boon2022mixed}
Wietse~M Boon and Jan~M Nordbotten.
\newblock Mixed-dimensional poromechanical models of fractured porous media.
\newblock \emph{Acta Mechanica}, pages 1--48, 2022.

\bibitem[Varela et~al.(2022)Varela, Ahmed, Keilegavlen, Nordbotten, and Radu]{varela2022aposteriori}
Jhabriel Varela, Elyes Ahmed, Eirik Keilegavlen, Jan~M Nordbotten, and Florin~A Radu.
\newblock A posteriori error estimates for hierarchical mixed-dimensional elliptic equations.
\newblock \emph{Journal of Numerical Mathematics}, 2022.

\bibitem[Coussy(2004)]{coussy2004poromechanics}
Olivier Coussy.
\newblock \emph{Poromechanics}.
\newblock John Wiley \& Sons, 2004.

\bibitem[Nikolaevskij(1990)]{nikolaevskij1990mechanics}
Viktor~N Nikolaevskij.
\newblock \emph{Mechanics of porous and fractured media}, volume~8.
\newblock World Scientific, 1990.

\bibitem[Garipov and Hui(2019)]{garipov2019thm}
TT~Garipov and MH~Hui.
\newblock Discrete fracture modeling approach for simulating coupled thermo-hydro-mechanical effects in fractured reservoirs.
\newblock \emph{International Journal of Rock Mechanics and Mining Sciences}, 122:\penalty0 104075, 2019.

\bibitem[Viswanath and Natarajan(1989)]{viswanath1989data}
DS~Viswanath and G~Natarajan.
\newblock Data book on the viscosity of liquids; hemisphere pub.
\newblock \emph{Corp.: New York, NY, USA}, 1989.

\bibitem[Barton et~al.(1985)Barton, Bandis, and Bakhtar]{barton1985strength}
Nick Barton, Strength Bandis, and K~Bakhtar.
\newblock Strength, deformation and conductivity coupling of rock joints.
\newblock In \emph{International journal of rock mechanics and mining sciences \& geomechanics abstracts}, volume~22, pages 121--140. Elsevier, 1985.

\bibitem[Geuzaine and Remacle(2009)]{geuzaine2009gmsh}
Christophe Geuzaine and Jean-Fran{\c{c}}ois Remacle.
\newblock Gmsh: A 3-d finite element mesh generator with built-in pre-and post-processing facilities.
\newblock \emph{International journal for numerical methods in engineering}, 79\penalty0 (11):\penalty0 1309--1331, 2009.

\bibitem[Aziz and Settari(1979)]{aziz1979petroleum}
Khalid Aziz and Antonin Settari.
\newblock \emph{Petroleum reservoir simulation}, volume 476.
\newblock 1979.

\bibitem[Aavatsmark(2002)]{aavatsmark2002introduction}
I.~Aavatsmark.
\newblock An introduction to multipoint flux approximations for quadrilateral grids.
\newblock \emph{Computational Geosciences}, 6:\penalty0 405--432, 2002.

\bibitem[Nordbotten(2016)]{nordbotten2016stable}
Jan~Martin Nordbotten.
\newblock Stable cell-centered finite volume discretization for biot equations.
\newblock \emph{SIAM Journal on Numerical Analysis}, 54\penalty0 (2):\penalty0 942--968, 2016.

\bibitem[Nordbotten and Keilegavlen(2021)]{nordbotten2021mpxa}
Jan~Martin Nordbotten and Eirik Keilegavlen.
\newblock An introduction to multi-point flux (mpfa) and stress (mpsa) finite volume methods for thermo-poroelasticity.
\newblock In \emph{Polyhedral Methods in Geosciences}, pages 119--158. Springer, 2021.

\bibitem[Hüeber and Wohlmuth(2005)]{Huber2005}
S.~Hüeber and B.~I. Wohlmuth.
\newblock A primal–dual active set strategy for non-linear multibody contact problems.
\newblock \emph{Computer Methods in Applied Mechanics and Engineering}, 194\penalty0 (27-29):\penalty0 3147--3166, 2005.
\newblock ISSN 00457825.
\newblock \doi{10.1016/j.cma.2004.08.006}.

\bibitem[Stefansson et~al.(2021{\natexlab{a}})Stefansson, Berre, and Keilegavlen]{Stefansson2021thm}
Ivar Stefansson, Inga Berre, and Eirik Keilegavlen.
\newblock A fully coupled numerical model of thermo-hydro-mechanical processes and fracture contact mechanics in porous media.
\newblock \emph{Computer Methods in Applied Mechanics and Engineering}, 386:\penalty0 114122, 2021{\natexlab{a}}.
\newblock ISSN 0045-7825.
\newblock \doi{10.1016/j.cma.2021.114122}.

\bibitem[H{\"u}eber and Wohlmuth(2005)]{hueber2005primal}
Stefan H{\"u}eber and Barbara~I Wohlmuth.
\newblock A primal--dual active set strategy for non-linear multibody contact problems.
\newblock \emph{Computer Methods in Applied Mechanics and Engineering}, 194\penalty0 (27-29):\penalty0 3147--3166, 2005.
\newblock \doi{10.1016/j.cma.2004.08.006}.

\bibitem[Berge et~al.(2020)Berge, Berre, Keilegavlen, Nordbotten, and Wohlmuth]{berge2020poroelastic_contact}
Runar~L Berge, Inga Berre, Eirik Keilegavlen, Jan~M Nordbotten, and Barbara Wohlmuth.
\newblock Finite volume discretization for poroelastic media with fractures modeled by contact mechanics.
\newblock \emph{International Journal for Numerical Methods in Engineering}, 121\penalty0 (4):\penalty0 644--663, 2020.

\bibitem[Paszke et~al.(2017)Paszke, Gross, Chintala, Chanan, Yang, DeVito, Lin, Desmaison, Antiga, and Lerer]{paszke2017automatic}
Adam Paszke, Sam Gross, Soumith Chintala, Gregory Chanan, Edward Yang, Zachary DeVito, Zeming Lin, Alban Desmaison, Luca Antiga, and Adam Lerer.
\newblock Automatic differentiation in pytorch.
\newblock 2017.

\bibitem[Abadi et~al.(2017)Abadi, Isard, and Murray]{abadi2017computational}
Mart{\'\i}n Abadi, Michael Isard, and Derek~G Murray.
\newblock A computational model for tensorflow: an introduction.
\newblock In \emph{Proceedings of the 1st acm sigplan international workshop on machine learning and programming languages}, pages 1--7, 2017.

\bibitem[Naumann(2011)]{naumann2011art}
Uwe Naumann.
\newblock \emph{The art of differentiating computer programs: an introduction to algorithmic differentiation}.
\newblock SIAM, 2011.

\bibitem[Stefansson and Varela(2023)]{stefansson2023source_code}
Ivar Stefansson and Jhabriel Varela.
\newblock Source code and results flexible and rigorous numerical modelling of multiphysics processes in fractured porous media using porepy, 2023.

\bibitem[Khan and Khan(2014)]{khan2014importance}
Mohd~Ehmer Khan and Farmeena Khan.
\newblock Importance of software testing in software development life cycle.
\newblock \emph{International Journal of Computer Science Issues (IJCSI)}, 11\penalty0 (2):\penalty0 120, 2014.

\bibitem[Burnstein(2006)]{burnstein2006practical}
Ilene Burnstein.
\newblock \emph{Practical software testing: a process-oriented approach}.
\newblock Springer Science \& Business Media, 2006.

\bibitem[Krekel et~al.(2023)Krekel, Oliveira, Pfannschmidt, Bruynooghe, Laugher, and Bruhin]{pytest71}
Holger Krekel, Bruno Oliveira, Ronny Pfannschmidt, Floris Bruynooghe, Brianna Laugher, and Florian Bruhin.
\newblock pytest 7.1.3, 2023.
\newblock URL \url{https://github.com/pytest-dev/pytest}.

\bibitem[Harris et~al.(2020)Harris, Millman, Van Der~Walt, Gommers, Virtanen, Cournapeau, Wieser, Taylor, Berg, Smith, et~al.]{harris2020array}
Charles~R Harris, K~Jarrod Millman, St{\'e}fan~J Van Der~Walt, Ralf Gommers, Pauli Virtanen, David Cournapeau, Eric Wieser, Julian Taylor, Sebastian Berg, Nathaniel~J Smith, et~al.
\newblock Array programming with {NumPy}.
\newblock \emph{Nature}, 585\penalty0 (7825):\penalty0 357--362, 2020.

\bibitem[Virtanen et~al.(2020)Virtanen, Gommers, Oliphant, Haberland, Reddy, Cournapeau, Burovski, Peterson, Weckesser, Bright, et~al.]{virtanen2020scipy}
Pauli Virtanen, Ralf Gommers, Travis~E Oliphant, Matt Haberland, Tyler Reddy, David Cournapeau, Evgeni Burovski, Pearu Peterson, Warren Weckesser, Jonathan Bright, et~al.
\newblock {SciPy} 1.0: fundamental algorithms for scientific computing in {P}ython.
\newblock \emph{Nature methods}, 17\penalty0 (3):\penalty0 261--272, 2020.

\bibitem[Kempf and Koch(2017)]{kempf2017system}
Dominic Kempf and Timo Koch.
\newblock System testing in scientific numerical software frameworks using the example of {DUNE}.
\newblock \emph{Archive of Numerical Software}, 5\penalty0 (1):\penalty0 151--168, 2017.

\bibitem[Roy(2005)]{roy2005review}
Christopher~J Roy.
\newblock Review of code and solution verification procedures for computational simulation.
\newblock \emph{Journal of Computational Physics}, 205\penalty0 (1):\penalty0 131--156, 2005.

\bibitem[Oberkampf and Roy(2010)]{oberkampf2010verification}
William~L Oberkampf and Christopher~J Roy.
\newblock \emph{Verification and validation in scientific computing}.
\newblock Cambridge University Press, 2010.

\bibitem[Stefansson and Keilegavlen(2023)]{Stefansson2023numerical}
Ivar Stefansson and Eirik Keilegavlen.
\newblock Numerical treatment of state-dependent permeability in multiphysics problems.
\newblock \emph{Water Resources Research}, page e2023WR034686, 2023.
\newblock \doi{10.1029/2023WR034686}.

\bibitem[Banshoya et~al.(2023)Banshoya, Berre, and Keilegavlen]{banshoya2023}
Ivar Banshoya, Inga Berre, and Keilegavlen.
\newblock Simulation of reactive transport in fractured porous media.
\newblock \emph{Transport in Porous Media}, 149:\penalty0 643--667, 2023.

\bibitem[Stefansson et~al.(2021{\natexlab{b}})Stefansson, Keilegavlen, Halld{\'o}rsd{\'o}ttir, and Berre]{stefansson2021propagation}
Ivar Stefansson, Eirik Keilegavlen, S{\ae}unn Halld{\'o}rsd{\'o}ttir, and Inga Berre.
\newblock Numerical modelling of convection-driven cooling, deformation and fracturing of thermo-poroelastic media.
\newblock \emph{Transport in Porous Media}, 140:\penalty0 371--394, 2021{\natexlab{b}}.

\bibitem[Dang et~al.(2022)Dang, Berre, and Keilegavlen]{dang2022poro_wing_cracks}
Hau~Trung Dang, Inga Berre, and Eirik Keilegavlen.
\newblock Two-level simulation of injection-induced fracture slip and wing-crack propagation in poroelastic media.
\newblock \emph{International Journal of Rock Mechanics and Mining Sciences}, 160:\penalty0 105248, 2022.
\newblock ISSN 1365-1609.
\newblock \doi{10.1016/j.ijrmms.2022.105248}.

\end{thebibliography}

\begin{appendices}
    \counterwithin{figure}{section}

    \section{Manufactured solution for system test\label{sec:exact_sol}}

    \begin{figure}[!h]
        \centering
        \includegraphics[width=0.5\textwidth]{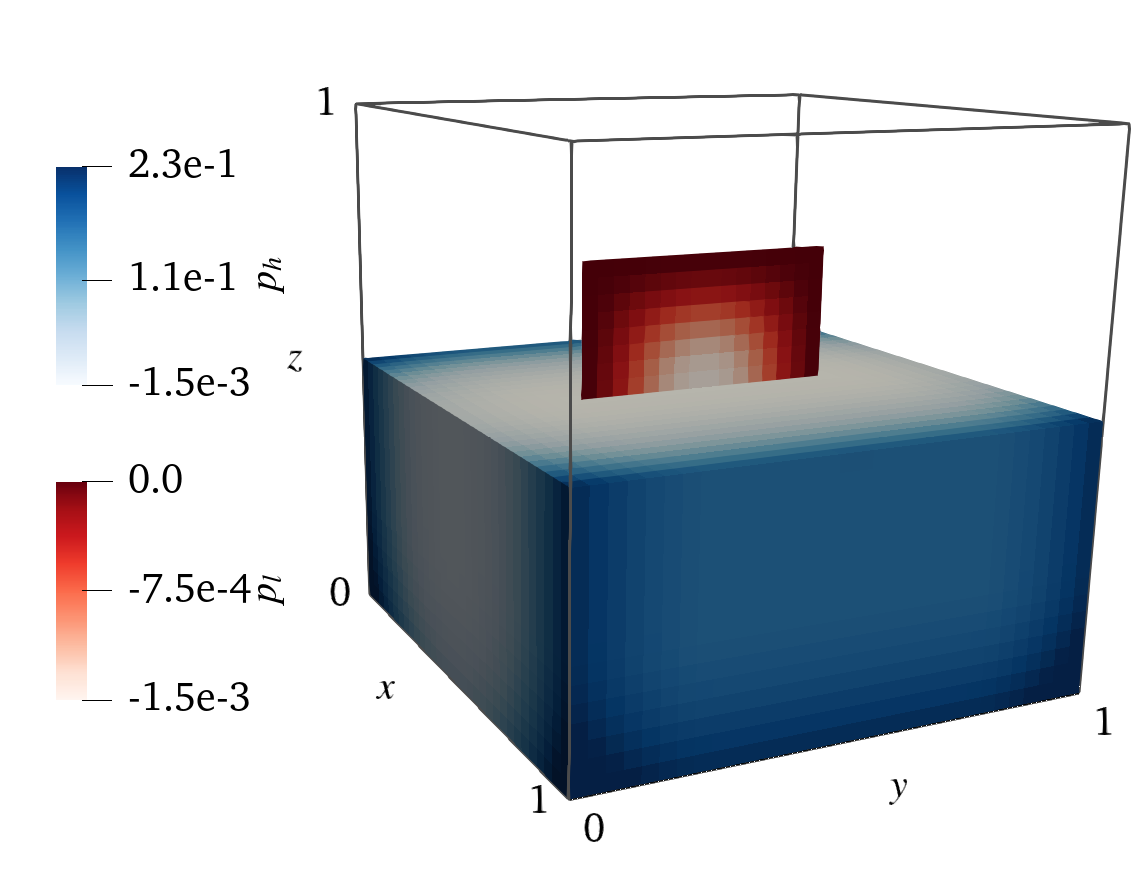}
        \includegraphics[width=0.48\textwidth]{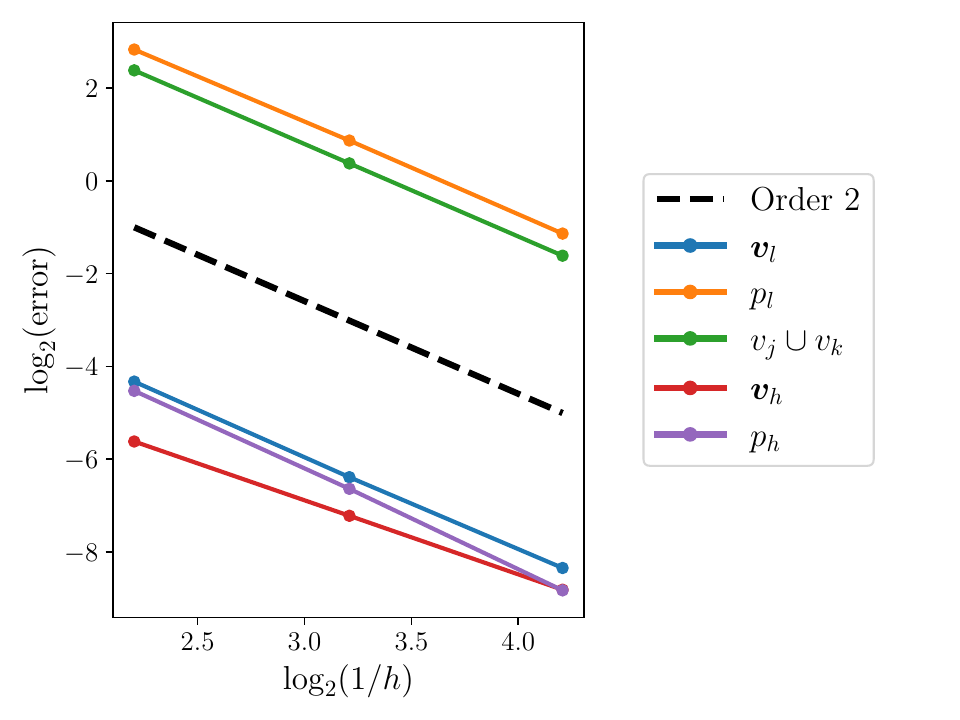}
        \caption{Analytical solution (left) and convergence analysis (right) for the compressible flow in fractured porous media using a Cartesian grid.
        The pressure solution is symmetric about the plane $x=\num{0.5}$, where the matrix grid is cut to expose the top half of the fracture.
        Rates for the interface fluxes are shown together, since \texttt{PorePy} treats two-sided interfaces as one single object.}
        \label{fig:compressible_flow}
    \end{figure}

    In this appendix, we provide the steps for deriving the manufactured solution of the system test presented in Section \ref{sec:system}. The solution follows closely the one presented in Appendix D.2 from \cite{varela2022aposteriori}, and can be seen as its generalisation to the compressible case. The final expressions, however, are considerably more involved and for that reason we do not include them explicitly. The interested reader can access all the expressions through the \texttt{ManuCompExactSolution3d} class from the \texttt{manu_flow_comp_3d_frac.py} module located in the PorePy sub-directory \texttt{porepy/tests/functional/setups}.

    As in \cite{varela2022aposteriori}, we consider a mixed-dimensional domain $Y = \Omega_l \sqcup \Omega_h = (0, 1)^3$ composed of a single vertical fracture
    \begin{equation}
        \Omega_l = \{\bm{x} \in Y: x = 0.5, 0.25 \leq y \leq 0.75, 0.25 \leq z \leq 0.75\},
    \end{equation}
    and a matrix $\Omega_h = Y \setminus \Omega_f$, which is divided into nine subregions, namely:
    \begin{equation}
        \Omega_h = \bigcup_{\alpha = 1}^9 \Omega_{h}^\alpha,
    \end{equation}
    where
    \begin{align}
        \Omega_h^1 &= \left\lbrace \bm{x}\in\Omega_h : 0.00 < y < 0.25,\,\, 0.00 < z < 0.25\right\rbrace, \\
        \Omega_h^2 &= \left\lbrace \bm{x}\in\Omega_h : 0.00 < y < 0.25,\,\, 0.25 \leq z < 0.75\right\rbrace, \\
        \Omega_h^3 &= \left\lbrace \bm{x}\in\Omega_h : 0.00 < y < 0.25,\,\, 0.75 \leq z < 1.00\right\rbrace, \\
        \Omega_h^4 &= \left\lbrace \bm{x}\in\Omega_h : 0.25 \leq y < 0.75,\,\, 0.00 < z < 0.25\right\rbrace, \\
        \Omega_h^5 &= \left\lbrace \bm{x}\in\Omega_h : 0.25 \leq y < 0.75,\,\, 0.25 \leq z < 0.75\right\rbrace, \\
        \Omega_h^6 &= \left\lbrace \bm{x}\in\Omega_h : 0.25 \leq y < 0.75,\,\, 0.75\leq z < 1.00\right\rbrace, \\
        \Omega_h^7 &= \left\lbrace \bm{x}\in\Omega_h : 0.75 \leq y < 1.00,\,\, 0.00 < z < 0.25\right\rbrace, \\
        \Omega_h^8 &= \left\lbrace \bm{x}\in\Omega_h : 0.75 \leq y < 1.00,\,\, 0.25 \leq z < 0.75\right\rbrace, \\
        \Omega_h^9 &= \left\lbrace \bm{x}\in\Omega_h : 0.75 \leq y < 1.00,\,\, 0.75 \leq z < 1.00\right\rbrace.
    \end{align}
    Moreover, we set the time interval of interest as $(0, \mathcal{T}) = (0, 1)$.

    The exact solutions are constructed with the help of the distance function $\delta(\bm{x})$, defined as the shortest length from any point in the matrix $\Omega_h$ to the fracture $\Omega_l$, and given by:
    \begin{equation}
        \delta(\bm{x}) :=
            \begin{cases}
            \sqrt{(x-0.5)^2 + (y-0.25)^2 + (z-0.25)^2}, & \Omega_h^1, \\
            \sqrt{(x-0.5)^2 + (y-0.25)^2}, & \Omega_h^2, \\
            \sqrt{(x-0.5)^2 + (y-0.25)^2 + (z-0.75)^2}, & \Omega_h^3, \\
            \sqrt{(x-0.5)^2 + (z-0.25)^2}, & \Omega_h^4, \\
            \sqrt{(x-0.5)^2}, & \Omega_h^5, \\
            \sqrt{(x-0.5)^2 + (z-0.75)^2}, & \Omega_h^6, \\
            \sqrt{(x-0.5)^2 + (y-0.75)^2 + (z-0.25)^2}, & \Omega_h^7, \\
            \sqrt{(x-0.5)^2 + (y-0.75)^2}, & \Omega_h^8, \\
            \sqrt{(x-0.5)^2 + (y-0.75)^2 + (z-0.75)^2}, & \Omega_h^9. \\
        \end{cases}
    \end{equation}

    We will also need the bubble function $\omega(\bm{x}) \in \Omega_h^5$:
    \begin{equation}
        \omega(\bm{x}) := 100 (y-0.25)^2 (y-0.75)^2 (z-0.25)^2 (z-0.75)^2.
    \end{equation}

    The manufactured solution is based on defining $p_h(\bm{x}, t)$  as a modified, smoother version of the distance function:
    \begin{equation}
        p_h(\bm{x}, t) :=
        t\begin{cases}
            \delta^{\xi+1}, & \left(\Omega_h \setminus \Omega_h^{5}\right) \times (0,\mathcal{T}) \\
            \delta^{\xi+1} + \omega\delta, & \Omega_h^5 \times (0,\mathcal{T}),
    \end{cases},\label{eq:appx_matrix_pressure}
    \end{equation}
    where $\xi \in \mathbb{R}_{>0}$ is a parameter that controls the regularity of the solution. Following \cite{varela2022aposteriori}, we employ $\xi = 1.5$, which offers sufficient smoothness while preserving non-trivial matrix/fracture coupling conditions.

    The density $\density{f}{h}(\bm{x},t)$ can now be obtained via Eq.~\eqref{eq:density}. For this particular test, we employ $\compressibility{f}{h} = 0.2$, $\pressure{0}{h} =0$ and $\density{f}{0}=1.0$. By setting $\permeability{}{h} / \viscosity{f}{h}= 1$, the Darcy's velocity $\fluidVelocity{}{h}(\bm{x}, t)$ can be obtained with the help of  \eqref{eq:darcy}. After setting $\porosity{}{h} = 0.1$, we have all the ingredients to compute the piece-wise time-dependent source term $\massSource{}{h}(\bm{x}, t)$ in the matrix via \eqref{eq:mass_conservation}.

    Due to continuity of normal mass fluxes \eqref{eq:internal_bc_mass} and noting that $\specificVolume{}{h}=\specificVolume{}{j}=\specificVolume{}{k}=1$, there holds:
    \begin{align}
         \density{f}{h}\fluidVelocity{}{h} \cdot \normal{}{h} &= \density{f}{j} \interfaceFluidFlux{}{j} = t \omega, \qquad \interface{}{j} \times (0, \mathcal{T}), \label{eq:appx_continuity_plus_side} \\
         \density{f}{h}\fluidVelocity{}{h} \cdot \normal{}{h} &= \density{f}{k} \interfaceFluidFlux{}{k} = t \omega, \qquad \interface{}{k} \times (0, \mathcal{T}),  \label{eq:appx_continuity_minus_side}
    \end{align}
    where $\interface{}{j}$ and $\interface{}{k}$ denote the interfaces coupling $\Omega_h$ and $\Omega_l$. Note that we have not included the projections operators in \eqref{eq:appx_continuity_plus_side} and \eqref{eq:appx_continuity_minus_side} since they trivially evaluate to identity matrices.

    It is straightforward to check that $p_h = 0$ on $\partial_j\Omega_h$ and $\partial_k\Omega_h$. Thus, setting $2\permeability{}{j} / \viscosity{f}{j} \aperture{}{l} =  2\permeability{}{k} / \viscosity{f}{k} \aperture{}{l}= 1$ in \eqref{eq:interfaceFluxes}, the fracture pressure is fixed and evaluates to the negative of the interface fluxes, i.e., $\pressure{}{l}(\bm{x}, t) = -t\omega$. With this, it is possible to obtain the fluid density in the fracture $\density{f}{l}(\bm{x}, t)$ via \eqref{eq:density}. As in the matrix, we set $\compressibility{f}{l} = 0.2$, $\density{0}{l}= 1$ and $\pressure{0}{l} = 0$. Assuming unit mobility $\permeability{}{l}/\viscosity{f}{l} = 1$, we can compute the tangential Darcy flux in the fracture $\fluidVelocity{}{l}(\bm{x}, t)$ via \eqref{eq:darcy}. By setting $\porosity{}{l} = 0.1$ and $\specificVolume{}{l}=1$, the time-dependent source term in the fracture $\massSource{}{l}(\bm{x},t)$ can now be obtained via \eqref{eq:mass_conservation}.

    To close the system of equations, we use zero initial conditions for all primary variables, impose Dirichlet boundary conditions satisfying \eqref{eq:appx_matrix_pressure} on the boundaries of $\Omega_h$ and no-flux at the tips of $\Omega_l$. Note that in the matrix, both pressure and density will change in every time step and must therefore be updated accordingly.

\end{appendices}

\end{document}